\magnification=\magstep1
\input amstex
\documentstyle{amsppt}
\catcode`\@=11 \loadmathfont{rsfs}
\def\mycal{\mathfont@\rsfs}
\csname rsfs \endcsname \catcode`\@=\active

\vsize=7.5in

\topmatter 
\title 
Coarse decomposition of II$_1$ factors    
\endtitle
\author SORIN POPA \endauthor

\rightheadtext{Coarse decomposition of factors}

\affil     {\it  University of California, Los Angeles} \endaffil

\address Math.Dept., UCLA, Los Angeles, CA 90095-1555, USA  \endaddress
\email popa\@math.ucla.edu \endemail

\thanks Supported in part by NSF Grant DMS-1700344 and the Takesaki Endowed Chair at UCLA \endthanks

\abstract  We prove that any separable II$_1$ factor $M$ admits a {\it coarse decomposition} over the hyperfinite 
II$_1$ factor $R$, i.e.,   
there exists an embedding $R\hookrightarrow M$ such that $L^2M\ominus L^2R$ is a multiple of the coarse Hilbert $R$-bimodule $L^2R \overline{\otimes} L^2R^{op}$.    
Equivalently, the von Neumann algebra generated by left and right multiplication by $R$ on $L^2M\ominus L^2R$ is isomorphic to $R\overline{\otimes}R^{op}$. 
Moreover, if $Q\subset M$ is an infinite index irreducible subfactor, then $R\hookrightarrow  M$ can be constructed to be coarse with respect to $Q$ as well.  
This implies existence of maximal abelian $^*$-subalgebras that are mixing, strongly malnormal, and with infinite multiplicity, in any given separable II$_1$ factor. 

\endabstract 

\endtopmatter

\document

\heading 1. Introduction \endheading

Structure and classification results in von Neumann  algebra framework depend heavily on the way an algebra  
decomposes as a bimodule over its subalgebras, especially over approximately finite dimensional (AFD) subalgebras,  
such as the hyperfinite II$_1$ factor $R=\overline{\otimes}_n (\Bbb M_2(\Bbb C), tr)_n$. For instance, if  a II$_1$ factor 
$M$ decomposes as a direct sum of  dichotomic classes of bimodules, 
compact (structured) and respectively coarse (mixing/random), then the ``tension'' between them enhances the analysis in $M$, leading to surprising rigidity phenomena (see e.g., [P06]). 

Along these lines, we prove in this paper a very general bimodule decomposition structure result. 
It shows that any separable II$_1$ factor $M$ contains a {\it coarse hyperfinite subfactor}, i.e.,  $M$ admits 
an embedding $R \hookrightarrow  M$ with the property that the Hilbert $R$-bimodule $_RL^2M_R$ 
decomposes as the direct sum between a copy of the trivial $R$-bimodule, $L^2R$, and a multiple of the coarse $R$-bimodule, $L^2R \overline{\otimes} L^2R^{op}$. 
Moreover, $R$ can be taken so that to satisfy several other ``constraints'', such as being contained 
in an irreducible subfactor $P\subset M$ and be almost orthogonal  
and coarse with respect to  a given subalgebra $Q\subset M$ satisfying $P\not\prec_M Q$  
(the {\it pair $R, Q$ is coarse} if $_RL^2M_Q$ is  a multiple of $L^2R \overline{\otimes} L^2Q^{op}$). 

The coarse subfactor $R\subset M$ is constructed as an inductive limit of dyadic matrix algebras, through an iterative technique also 
used in [P81a], [P16], [P17]. But while in all previous work the resulting bimodule structure  $_RL^2M_R$ remained 
``blind'', the big novelty here is that we are able to construct embeddings $R \hookrightarrow M$ with complete control of the 
bimodule decomposition  at the end of the iterative process. 

Coarseness  is in some sense the ``most random'' position a subalgebra  may have in the ambient II$_1$ factor. It automatically entails mixingness, 
which in turn implies strong malnormality, a property that's in dichotomy with the weak quasi regularity of the subalgebra. 
Altogether,  our main result shows the following: 

\proclaim{1.1. Theorem} Any separable  $\text{\rm II}_1$ factor  $M$ 
contains a hyperfinite factor $R\subset M$ 
that's coarse in $M$ $($and thus also mixing and strongly malnormal in $M)$. 
Moreover, given any irreducible subfactor $P\subset M$, any von 
Neumann subalgebra $Q\subset M$ satisfying $P\not\prec_M Q$ and any $\varepsilon >0$, the coarse  subfactor $R\subset M$ 
can be constructed so that to be contained in $P$,  be coarse with respect to $Q$ and satisfy $R \perp_\varepsilon Q$. 
\endproclaim

The condition $P\not\prec_M Q$ for two subalgebras of the II$_1$ factor $M$ is in the sense of 
(Definition 2.4 in [P03]) and means that there exists no non-zero {\it intertwiner} from $P$ to $Q$ (i.e., $x\in M$ with   
$\text{\rm dim}(L^2(PxQ)_Q)<\infty$). This automatically implies that $Q$ has 
uniform infinite index in $M$, i.e.,  
given any non-zero projection $p\in Q'\cap M$, 
the  [PiP84]-index of the inclusion $Qp\subset pMp$ is infinite. When $Q$ is an irreducible subfactor of $M$, 
it amounts to $Q\subset M$ having infinite Jones index [J83].

A subalgebra $B\subset M$ is {\it mixing} if the action Ad$\Cal U(B) \curvearrowright M$ is mixing relative to $B$ 
in the sense of (2.9 in [P05]), equivalently if $\lim_u \|E_B(xuy)\|_2 =0$, for all $x, y\in M\ominus B$, where the limit is over $u\in \Cal U(B)$ tending weakly to $0$. 
The subalgebra $B\subset M$ is {\it strongly malnormal} if 
any $x\in M$ that's a {\it weak intertwiner} for $B$, i.e., satisfies $\text{\rm dim}(L^2(A_0xB)_B)<\infty$, for some diffuse $A_0\subset B$,  
lies in $B$  (see [P04], [IPeP05], [PeT07], [GP14] for variations of this property for subgroups and subalgebras).

One can show that if $R\subset M$ is coarse then $R\subset M$ is mixing (see 2.6.3 $(a)$ below). In turn, by (3.1 in [P03]), 
the mixing property implies very strong 
absorption properties for $R\subset M$, meaning that $R$ is strongly malnormal in the above sense.  
In particular, any maximal abelian $^*$-subalgebra (abbreviated hereafter as {\it MASA}) $A$ 
of $R$ is a MASA in $M$, with all its weak intertwiners  contained in $R$. 
Moreover, if a MASA $A\subset R$ is  coarse in $R$ with 
$L^2R\ominus L^2A \simeq (L^2A \overline{\otimes} L^2A)^{\oplus\infty}$ as Hilbert $A$-bimodules, then $A$ is coarse in $M$ as well, 
and if $A$ has infinite multiplicity in $R$, then so it does in $M$ (i.e.,  $(A \vee A^{op})'$ is of type I$_\infty$ on $L^2M\ominus L^2A$; 
the type of $(A \vee A^{op})'\cap \Cal B(L^2M\ominus L^2A)$ is sometimes called the {\it Pukanszky invariant} of $A\subset M$, see [Pu60] or [P16]). 
Thus, if one represents  a coarse hyperfinite  II$_1$ subfactor $R\subset M$ 
as the II$_1$ factor of the lamp-lighter group $R=L(\Bbb Z/2\Bbb Z \wr \Bbb Z)$, 
then $A=L(\Bbb Z)$ follows  coarse,  strongly malnormal, with infinite multiplicity  in $M$. So we have:

\proclaim{1.2. Corollary}  Any separable $\text{\rm II}_1$ factor $M$ 
contains a MASA $A\subset M$ that's coarse in $M$ $($so in particular strongly malnormal and mixing$)$ and which in addition can be taken to have infinite multiplicity in $M$. 
Moreover, given any irreducible subfactor $P\subset M$, any von Neumann subalgebra $Q\subset M$ 
such that $P\not\prec_M Q$ and any $\varepsilon >0$, the coarse MASA $A\subset M$ can be constructed inside $P$, coarse to $Q$, and satisfying $A\perp_\varepsilon Q$. 
\endproclaim

The problems of whether any separable II$_1$ factor contains malnormal MASAs and MASAs with infinite multiplicity,  
both of which are strengthening of singularity, have been open  for some time. The second of these problems  
is implicit in ([Pu60], [P81c]) and has been largely circulated as a folklore problem.  
Strengthened singularity such as mixingness and malnormality is discussed in (Section 5.1 in [P13a] and Section 5.3 in [P16]). 
It is somewhat surprising that any II$_1$ factor contains a singular MASA having all of these properties.

Theorem 1.1 also implies that any Cartan subalgebra $D$ of a coarse  subfactor $R\subset M$ is maximal abelian 
in $M$, with its normalizer in $M$ generating $R$, thus showing   
existence of semiregular MASAs $D\subset M$ whose normalizing algebra is a hyperfinite II$_1$ factor $R\subset M$  
which in addition can be taken  $\varepsilon$-orthogonal and coarse with respect  to some given irreducible infinite index subfactor $Q\subset M$. 
On the other hand, since any tracial AFD algebra can be embedded into $R$ and any countable amenable group $G$  
gives rise to a tracial AFD  von Neumann algebra (by [C76]), this also shows existence of copies of the left regular representation of $G$ that 
are $\perp_\varepsilon Q$. 

\proclaim{1.3. Corollary} Let  $M$ be a separable $\text{\rm II}_1$ factor, $P\subset M$ an irreducible subfactor and $Q\subset M$ 
a von Neumann subalgebra such that $P\not\prec_M Q$.  Let also $\varepsilon >0$. 

$1^\circ$ There exists a semiregular MASA $D$ of $M$  that's contained in $P$, whose normalizer $\Cal N_M(A)$ lies in $P$ and generates a hyperfinite factor $R$ 
satisfying $R \perp_\varepsilon Q$. 

$2^\circ$ If $G$ is a countable amenable group, then there exists a copy $\{u_g\}_{g\in G}\subset P$ of the left regular representation of $G$ such that 
$\|E_Q(u_g)\|_2 \leq \varepsilon$, $\forall g\in G\setminus \{e\}$. 
\endproclaim

The $\varepsilon$-orthogonality   between subalgebras in these statements is with respect to the Hilbert structure given by the (unique) trace state $\tau$ on the ambient factor $M$. 
Thus, $B \perp_\varepsilon Q$ 
for subalgebras $B, Q\subset M$ 
means that $\|E_Q(b)\|_2 \leq \varepsilon \|b\|_2$ for all $b\in B\ominus \Bbb C\overset{def}\to{=} \{b\in B \mid \tau(b)=0\}$, where as usual $E_Q$ denotes the 
trace preserving expectation onto $Q$. 
We will in fact prove $\varepsilon$-perpendicularity in a stronger sense, with the hyperfinite II$_1$ factor $R\subset P$ in Theorem 1.1   
constructed so that  $\|E_Q(b)\|_q \leq \varepsilon \|b\|_q$ for all $b\in R\ominus \Bbb C1$ and all $1\leq q \leq 2$ (see Theorem 4.2).  This 
is equivalent to the condition $\|E_R(x)\|_p\leq \varepsilon \|x\|_p$, $\forall x\in Q\ominus \Bbb C1$, $\forall 2\leq p \leq \infty$, so in particular $R$ satisfies 
$\|E_R(x)\|\leq \varepsilon \|x\|$, $\forall x\in Q\ominus \Bbb C1$.

As we mentioned before, the construction of  the coarse hyperfinite II$_1$ subfactor $R\subset M$ in Theorem 1.1 
uses the iterative strategy from ([P81a], [P81c], [P16]). Thus, $R$ is obtained as an inductive limit of I$_{2^n}$-subfactors $B_n\simeq 
\Bbb M_2(\Bbb C)^{\otimes n}$ inside $P$, so that at each ``next step'' $n+1$ the ``newly added'' algebra $B_n'\cap B_{n+1}\simeq \Bbb M_2(\Bbb C)$ becomes almost 2-independent to $M\ominus B_n$, while 
remaining almost orthogonal to $Q$. 
We do this by using incremental patching technique ([P92] ,[P13a], [P13b], [P17]), which we combine with a key new technique of controlling the states implemented on 
$B_{n+1} \vee B_{n+1}^{op}$ by finitely many given vectors in $L^2M$. If the asymptotic $2$-independence is made rapidly enough, and with more and more elements in $M$ being taken into account, 
then the resulting $R$ is so that $R\otimes_{C^*} R^{op} \subset \Cal B(L^2M \ominus L^2R)$ extends to a normal representation of $R\overline{\otimes} R^{op}$ 
on $L^2M\ominus L^2R$. This amounts  to $R\subset M$ being coarse.

In Section 2 of the paper we recall some definitions 
and prove some preliminary technical results. In Section 3 we prove the main technical Lemma 3.1 
needed in the proof of Theorem 1.1. Section 4 contains the actual proof of this theorem, in fact a stronger result is shown (see Theorem 4.2).  We end that section with 
several remarks,  including a definition of  {\it multiple mixing}  for an inclusion $B\subset M$, for which a 
result similar to Theorem 1.1 holds true. 

The technique for constructing AFD-embeddings with control of bimodule structure introduced in this paper 
opens up new perspectives for revisiting some well known unanswered questions in II$_1$ factors. We  dwell on this in Section 5, where we discuss two specific problems. 
Thus, we first conjecture that if any amplification of $M$ can be generated by two selfadjoint elements, then 
one can construct hyperfinite factors $R_0, R_1 \subset M$ recursively so that $R_0\vee R_1^{op}$ is purely infinite (i.e., has no finite summand). 
If true, then due to results in ([V88], [R92],  [V96], [GeP98]), this fact would imply 
that $M=L(\Bbb F_\infty)$ cannot be finitely generated, and that all free group factors are non-isomorphic. 
We then discuss the possibility that maximal AFD  
subalgebras of $L(\Bbb F_n)$ are all coarse (see also 4.3 in [H15]), and mutually coarse one to another, unless unitary conjugate.  
We prove that this property would imply that any maximal AFD is mixing (resp. strongly malnormal), which in turn would imply 
that any two maximal amenable subalgebras with diffuse intersection must coincide, a property that has been predicted in ([PeT07]).

\heading 2. Some preliminaries  \endheading

We use in this paper the same notations as in ([P13b], [P16], [P17])  and refer to [AP17]  
for basics in II$_1$ factors theory. One notation that's frequently used is $(\Cal B)_r$ for the ball of radius $r>0$ of a given Banach space 
$\Cal B$ (the space and its norm being clear from the context).

Recall that a tracial von Neumann algebra $(M, \tau)$ is said to be {\it separable} if 
$M$ is separable in the Hilbert norm $\| x \|_2 = \tau(x^*x)^{1/2}$ implemented by the (fixed) 
normal faithful trace state $\tau$. For tracial von Neumann algebras, this condition is equivalent to $M$ 
being countably generated (see e.g., [AP17]).

\vskip.05in 
\noindent
{\bf 2.1. The ultrapower framework.} We will often use the ultrapower formalism: If $M$ is a II$_1$ factor and $\omega$ is a free 
(or non-principal) ultrafilter $\omega$ on $\Bbb N$, then $M^\omega$ denotes its $\omega$-{\it ultrapower}  
II$_1$ factor (see e.g.  [C76], [AP17], or 1.6 in [P13b] for complete definitions). Thus, $M^\omega$ 
is endowed with the ultrapower trace $\tau((x_n)_n)=\underset{n\rightarrow\omega}\to{\lim} \tau(x_n)$, $\forall \  (x_n)_n \in M^\omega$. 

Also, if $N$ is a von Neumann subalgebra of the II$_1$ factor $M$, then we denote by $e_N$ the orthogonal projection 
of $L^2M$ onto $L^2N$. Thus, if $x\in M$ is viewed as the operator of left multiplication 
on $L^2M$, then $e_N x e_N = E_N(x)e_N$, forall $x\in M$, implying that sp$Me_N M$ 
is a $^*$-subalgebra in $\Cal B(L^2M)$. Following Jones notations and terminology 
from ([J83]),  we denote by $\langle M, e_N \rangle \subset \Cal B(L^2M)$ 
the {\it basic construction} algebra 
$\overline{\text{\rm sp}M e_N M}^w=(J_M N J_M)'$, where $J_M$ is the canonical conjugation on $L^2M$, 
$J(\xi)=\xi^*$, $\forall \xi \in L^2M$. We have $e_N \langle M, e_N \rangle e_N=Ne_N$ and 
the semi-finite 
von Neumann algebra $\langle M, e_N \rangle$ is endowed with 
the canonical normal faithful semi-finite trace $Tr=Tr_{\langle M, e_N \rangle}$, 
satisfying the condition $Tr(xey)=\tau(xy)$, for $x, y\in M$.

If $N\subset M$ is a von Neumann subalgebra and we consider the 
corresponding inclusion of ultrapower algebras $N^\omega \subset M^\omega$, 
then it is useful to keep in mind that the canonical trace $Tr$ of elements in the basic construction algebra $\langle M^\omega, e_{N^\omega} \rangle$ 
is obtained as a limit of the trace $Tr$ of elements in $\langle M, e_N \rangle$: 

\proclaim{2.1.1. Lemma} If $x=(x_n)_n, y=(y_n)_n \in M^\omega$ then $e_{N^\omega}xe_{N^\omega}=E_{N^\omega}(x)e_{N^\omega}=
(E_N(x_n))_ne_{N^\omega}$ and 
$$
Tr(xe_{N^\omega}y)=\lim_{n\rightarrow \omega} Tr(x_ne_Ny_n)=\lim_{n\rightarrow \omega} \tau(x_ny_n).
$$ 

\endproclaim
\noindent
{\it Proof}. Immediate by the definitions. 
\hfill 
$\square$ 

\vskip.05in 
\noindent
{\bf 2.2. Intertwining subalgebras}. Recall from (Section 2 in [P03]) that if $Q, P$ are von Neumann subalgebras 
of a tracial von Neumann algebra $M$, then the notation $Q\prec_M P$ means that 
there exists a non-zero $x\in M$ such that the Hilbert $Q-P$ bimodule $L^2(QxP)$ has finite dimension 
as a right $P$-module.  Following ([P16]; cf. also [P05a]), such $x$ will be called {\it intertwiners} from $Q$ to $P$ and   
we denote by $\Cal I_M(Q, P)$ the space of all such $x$, calling it  
the {\it intertwining space}  from $Q$ to $P$. Thus, $Q\prec_M P$ means that $\Cal I_M(Q, P)\neq 0$. 

In turn, if $\Cal I_M(Q, P)=0$ then we write $Q\not\prec_M P$. By (2.1-2.4 in [P03]; see also Section 1.3 in [P16]) this is equivalent to:  
$\forall F \subset M$ finite, $\forall \varepsilon >0$, $\exists u\in \Cal U(P)$ such that $\|E_Q(xuy)\|_2 \leq \varepsilon$, $\forall x, y \in F$. This last condition readily implies 
that if $P\not\prec_M Q$ then $P^\omega\not\prec_{M^\omega} Q^\omega$ (see e.g., 2.1 in [P17]; N.B. the converse holds true as well). 

It is trivial to see, by using the definitions, 
that if $P\not\prec_M Q$ and $B_0\subset P$ is a finite dimensional $^*$-subalgebra, 
then $(B_0'\cap P) \not\prec_M Q$ as well. When passing to ultra powers, this entails: 

\proclaim{2.2.1. Lemma} With the above notations, if $B\subset P^\omega$ is a separable AFD von Neumann subalgebra, 
then $P\not\prec_M Q$ implies $(B'\cap P^\omega)\not\prec_{M^\omega} Q^\omega$.  
\endproclaim
\noindent
{\it Proof}. Let $B_n\subset R$ be an increasing sequence of finite dimensional subalgebras that generates $B$. 
Let $F\subset (M^\omega)_1$ be a finite set and $\varepsilon >0$. 

Since $(B_n'\cap P)^\omega=B_n'\cap P^\omega \not\prec_{M^\omega} Q^\omega$, there exists a unitary element $u_n=(u_{n,k})_k\in B_n'\cap P^\omega$,  
with $u_{n,k}\in \Cal U(B_n'\cap P)$, $\forall k$, such that $\|E_{Q^\omega}(xu_ny)\|_2\leq \varepsilon/2$, for all $x, y\in F$. Thus, if $x=(x_k)_k, y=(y_k)_k$ 
with $x_k, y_k\in (M)_1$, then 
$$
\lim_{k\rightarrow \omega} \|E_Q(x_k u_{n,k}y_k\|_2 \leq \varepsilon/2 < \varepsilon. 
$$

Denote by $V_n$ the set of all $k\in \Bbb N$ such that   $ \|E_Q(x_k u_{n,k}y_k\|_2 < \varepsilon$ for all $x, y\in F$. 
Note that $V_n$ corresponds to an open closed neighborhood of $\omega$ in $\Omega$, 
under the identification $\ell^\infty\Bbb N = C(\Omega)$. Let now $W_n\subset \Bbb N$, $n\geq 0$, be defined recursively as follows: 
$W_0=\Bbb N$ and $W_{n+1}=W_n \cap V_{n+1}\cap \{k \in \Bbb N \mid k > \min W_n\}$. 
Note that, with the same identification as before, $W_n$ is a strictly decreasing sequence of neighborhoods of $\omega$ in $\Omega$. 

Define $u=(u'_m)_m$, where $u'_k=u_{m, k}$ for $k\in W_{m-1}\setminus W_m$. Then the above conditions show that $u$ 
is a unitary element in $\cap_n B_n'\cap P^\omega = B'\cap P^\omega$ which satisfies $\|E_{Q^\omega}(xuy)\|_2\leq \varepsilon$, $\forall x, y\in F$. 
\hfill 
$\square$ 

\vskip.05in 
\noindent
{\bf 2.3. Approximate $2$-independence of subalgebras}. Recall from ([P13a], [P13b])  that if $B_0, B\subset M$ are von Neumann algebras and 
$E\subset M\ominus B_0\overset{def}\to{=}\{x\in M \mid E_{B_0}(x)=0\}$ is a subset, then $B$ is $n$-independent 
to $E$ relative to $B_0$, if the expectation on $B_0$ of any word with alternating letters from $E, B\ominus \Bbb C$ and length at most 2$n$ 
(so at most $n$ alternations) is equal to $0$. 

\proclaim{2.3. Lemma}  Assume $P\subset M$ is an irreducible inclusion of $\text{\rm II}_1$ factors and $P_0\subset P^\omega$ a finite dimensional factor. 
Given any finite set $E\subset M^\omega \ominus P_0$, there exists a diffuse abelian subalgebra $A\subset P_0'\cap P^\omega$ such 
that $A$ is free independent to $E$ relative to $P_0$. In particular, $A$ is $2$-independent to $E$ relative to $P_0$ and thus, for any $1\geq c >0$ 
there exists a projection $q\in P_0'\cap P^\omega$ of trace $c$ such that $E_{P_0}(qz)=0$ and $E_{P_0}(qzqz^*)=\tau(q)^2E_{P_0}(zz^*)=c^2E_{P_0}(zz^*)$, $\forall z\in E$. 
\endproclaim
\noindent
{\it Proof}. This is just a particular case of (Lemma 1.4 in [P92]) or (Theorem 4.3 in [P13b]).  
\hfill 
$\square$

\vskip.05in 
\noindent
{\bf 2.4. Almost $L^p$-orthogonality of subalgebras}.  Recall now that if $y\in M$ and $1\leq p < \infty$, 
then one denotes $\|y\|_p=\tau(|y|^p)^{1/p}$. For a fixed $y$, the $L^p$-norms $\|y\|_p$ are increasing in $p$,  with the limit $\underset{p\rightarrow \infty}\to{\lim}\|y\|_p$ 
equal to the operator norm $\|y\|$, which we also view as $\|y\|_\infty$. The completion $L^p(M)$ of $M$ 
in the norm $\| \ \|_p$ identifies naturally with the  space of densely defined closed operators $Y$ on $L^2M$ that are affiliated with $M$ 
and have the property that $|Y|$ has spectral decomposition $|Y|=\int \lambda \text{\rm d}e_\lambda$ satisfying 
$\int \lambda^p \text{\rm d}\tau(e_\lambda) < \infty$. 

It is well known that if $1\leq p <\infty$ then $(L^pM)^*\simeq L^qM$,  where $q=\frac{p}{p-1}$ (with the usual convention 
$1/0=\infty$), the duality being given by $(\xi, \zeta) \mapsto \tau(\zeta \xi)$ 
for $\xi \in L^pM$, $\zeta\in L^qM$, viewed as operators affiliated with $M$. This also shows that 
if $y\in M$ and $1\leq p, q \leq \infty$ with $\frac{1}{p}+\frac{1}{q}=1$, then $\|y\|_p=\sup \{|\tau(yz)| \mid z\in (L^qM)_1 \}$. 

It is also useful to recall that if $x\in M\simeq \Bbb M_n(\Bbb C)$  and $1\leq p \leq  p' \leq \infty$ 
then $\|x\|_{p'} \leq  n^{\frac{1}{p}-\frac{1}{p'}}\|x\|_p\leq n \|x\|_p$. 

If $Q\subset M$ is a von Neumann algebra, then  $\tau(xy)=\tau(xE_Q(y))$, for all $x\in Q$, $y\in M$. So the above formula 
for calculating $\| \ \|_p$ shows that $\|E_Q(y)\|_p \leq \|y\|_p$, $\forall 1\leq p \leq \infty$. 

\vskip.05in
\noindent
{\it 2.4.1. Notation}. Let $B, Q$ be von Neumann subalgebras of the II$_1$ factor $M$ and $1\leq p \leq \infty$. 
For each $y\in M$, we denote $\text{\bf c}_p(yBy^*,Q)=\sup \{\|E_{Q}(yby^*)\|_p/\|b\|_p \mid b\in B\ominus \Bbb C1, b\neq 0\}$. Note that by the above remarks we have  
$\text{\bf c}_p(yBy^*,Q)=\sup \{|\tau(yby^*x)| \mid b\in B\ominus \Bbb C, \|b\|_p\leq 1, x\in Q, \|x\|_q\leq 1\}$,  where $q=\frac{p}{p-1}$. We'll also use the related 
constant $\text{\bf c}'_p(yBy^*,Q)=\sup \{|\tau(yby^*x)| \mid b\in B\ominus \Bbb C, \|b\|_p\leq 1, x\in Q\ominus \Bbb C, \|x\|_q\leq 1\}$, 
which clearly satisfies $\text{\bf c}'_p(yBy^*,Q)= \text{\bf c}'_q(y^*Qy, B)$. 

\vskip.1in

\proclaim{2.4.2. Lemma}  With the above notations, we have 

$1^\circ$ $\text{\bf c}'_p(B,Q)\leq \text{\bf c}_p(B,Q) \leq 2\text{\bf c}'_p(B,Q) = 2 \text{\bf c}'_q(Q,B)\leq 2 \text{\bf c}_q(Q,B)$. 

$2^\circ$   If $B\simeq \Bbb M_{n}(\Bbb C)$ and $1\leq p \leq p' \leq \infty$, then $\text{\bf c}_{p'}(Q, B) \leq n^{\frac{1}{p}-\frac{1}{p'}} \text{\bf c}_p(Q, B)$. 

$3^\circ$ If $B_n\subset M$ is an increasing sequence of von Neumann algebras and $B=\overline{\cup_n B_n}^w$, then $\underset{n\rightarrow \infty}\to{\lim} \text{\bf c}_p(B_n,Q)
=\text{\bf c}_p(B,Q)$. 
\endproclaim
\noindent
{\it Proof}. $1^\circ$ If $b\in B\ominus \Bbb C$ and $x\in Q$, then  $\tau(bx)=\tau(b(x- \tau(x)1))$ and $\|x-\tau(x)\|_q \leq 2\|x\|_q$. 
Thus, if $x\in (L^qQ)_1$ then $\|x-\tau(x)1\|_q \leq 2\|x\|_q \leq 2$ and so we have 
$$
\text{\bf c}_p(B,Q)=\sup\{|\tau(bx)| \mid b\in (L^pB \ominus \Bbb C)_1, x\in (L^qQ)_1\} 
$$
$$
\leq \sup\{|\tau(by)| \mid b\in (L^pB \ominus \Bbb C)_1, y\in (L^qQ \ominus \Bbb C)_2\} 
$$
$$
= 2 \sup\{|\tau(by)| \mid b\in (L^pB \ominus \Bbb C)_1, y\in (L^qQ \ominus \Bbb C)_1\} = 2\text{\bf c}'_p(B,Q)
$$
$$
= 2 \text{\bf c}'_q(Q,B)= 2 \sup\{|\tau(yb)| \mid  y\in (L^qQ \ominus \Bbb C)_1, b\in (L^pB \ominus \Bbb C)_1\}
$$
$$
\leq 2 \sup\{|\tau(yb)| \mid  y\in (L^qQ \ominus \Bbb C)_1, b\in (L^pB)_1\}=2\text{\bf c}_q(Q,B).
$$

$2^\circ$ Let $q=\frac{p}{p-1}$ and $q'=\frac{p'}{p'-1}$ and note that $1\leq q'\leq q \leq \infty$. Since the unit ball of $L^{p'}Q$ is included in the unit ball of $L^pQ$,  the uint ball of $L^{q'}B$ 
is included in the ball of radius $n^{\frac{1}{q'}-\frac{1}{q}}$ of $L^{q}B$ and we have $\frac{1}{q'}-\frac{1}{q}=\frac{1}{p}-\frac{1}{p'}$, it follows that 
$$
\text{\bf c}_{p'}(Q, B)= \sup\{|\tau(yb)| \mid  y\in (L^{p'}Q \ominus \Bbb C)_1, b\in (L^{q'}B)_1\} 
$$
$$
\leq  n^{\frac{1}{p}-\frac{1}{p'}} \sup\{|\tau(yb)| \mid  y\in (L^pQ \ominus \Bbb C)_1, b\in (L^qB)_1\}=n^{\frac{1}{p}-\frac{1}{p'}} \text{\bf c}_p(Q, B). 
$$

$3^\circ$ is  straightforward and we leave its proof as an exercise. 
\hfill 
$\square$ 
\vskip.05in 

\noindent
{\bf 2.5. Coarse subalgebras and coarse pairs}. We consider here a new property for 
a subalgebra of a II$_1$ factor, as well as for pairs of subalgebras. 

We recall in this respect that if $B, B_0$ are tracial von Neumann algebras 
then $L^2B\overline{\otimes} L^2B_0^{op}\simeq L^2(B\overline{\otimes}B_0^{op})$, with its $B-B_0$ bimodule structure given by $b(\xi \otimes \eta^{op})=(b\xi) \otimes \eta^{op}$ 
and $(\xi \otimes \eta^{op})b_0=\xi\otimes (\eta b_0)^{op}$, $\forall \xi \in L^2B, \eta \in L^2B_0$, $b\in B, b_0\in B_0$, 
is called the {\it coarse Hilbert} $B-B_0$ {\it bimodule}. 

If $p$ is a projection in $B \overline{\otimes} B_0^{op}$, and we denote by $\rho(p)$ the right multiplication by $p$ on $L^2(B\overline{\otimes} B_0^{op})$, 
then $\rho(p)(L^2(B\overline{\otimes} B_0^{op}))$ is still a Hilbert $B-B_0$ bimodule. We say that $_B\Cal H_{B_0}$ is a {\it multiple 
of the coarse} $B-B_0$ bimodule, if $\Cal H$ is a direct sum of such bimodules, $\Cal H=\oplus_i \rho(p_i)(L^2(B\overline{\otimes}B_0^{op}))$. 
Note that this is equivalent to the fact that 
the von Neumann algebra $B\vee B_0^{op}\subset \Cal B(\Cal H)$, generated by the operators of left multiplication by $B$ and right 
multiplication by $B_0$ on $\Cal H$, extends to  a (normal) representation of the von Neumann algebra $B\overline{\otimes} B_0^{op}$.

\vskip.05in
\noindent
{\it 2.5.1. Definition}. A (proper) diffuse von Neumann subalgebra $B$ of a tracial von Neumann algebra $M$ is a {\it coarse subalgebra} of $M$  if 
the Hilbert $B$-bimodule $L^2M \ominus L^2B$ is a multiple of the coarse $B$-bimodule $L^2B\overline{\otimes} L^2B^{op}$, or equivalently 
if  the  von Neumann algebra generated by $B$ and its mirror image 
$B^{op}=J_MBJ_M$ on $L^2M \ominus L^2B$ generate a normal representation of the von Neumann algebra $B \overline{\otimes} B^{op}$. 
We then also say that $M$ has a {\it coarse decomposition over $B$}. 
If  $B$ is AFD (e.g., $B\simeq R$), then an alternative terminology is that $B \hookrightarrow M$ is a 
{\it coarse embedding} of $B$ into $M$. 

Note that in the above definition of $B\subset M$ being coarse we have not assumed 
the faithfulness of the normal representation of $B\overline{\otimes} B^{op}$ in $\Cal B(L^2M\ominus L^2B)$. However, faithfulness  is 
automatic if $M$ is a II$_1$ factor. In fact, it fails only if $B, M$ have a common central projection on which they coincide  (see Proposition 2.5.4 ). 

\vskip.05in
\noindent
{\it 2.5.2. Definition}. Let $B, B_0\subset M$ be diffuse weakly closed $^*$-subalgebras of the tracial von Neumann algebra $M$, 
with supports $q=1_B$, $q_0=1_{B_0}$. We say that 
$B$ is {\it coarse with respect to} $B_0$ in $M$ if the Hilbert-bimodule $_{B}(qL^2Mq_0)_{B_0}$ is isomorphic to a multiple of the 
coarse $B-B_0$ bimodule $L^2B\overline{\otimes} L^2B_0^{op}$, or equivalently, if the von Neumann 
algebra generated by $B$ and $B_0^{op}=J_MB_0J_M$ in $\Cal B(qL^2Mq_0)$ is a normal representation of $B\overline{\otimes} B_0^{op}$. We then also say 
that $B, B_0$  is a {\it coarse pair} in $M$. 
(N.B. As with the definition of coarse subalgebra, we are not assuming faithfulness of $B\overline{\otimes} B_0^{op}\subset \Cal B(L^2M)$).

\vskip.05in
\noindent
{\it 2.5.3. Examples} $1^\circ$ If the tracial von Neumann algebra $M$ arises from an infinite  group $\Gamma$ 
and $H\subset \Gamma$ is an infinite subgroup, then the inclusion $B=L(H)\subset L(\Gamma)=M$ 
is coarse if and only if for any $g\in \Gamma \setminus H$ one has $gHg^{-1}\cap H = \{e\}$. Also, if $H_0\subset \Gamma$ is another group, 
then $L(H)$ and $L(H_0)$ is a coarse pair if and only if $gHg^{-1}\cap H_0=\{e\}$, for all $g\in \Gamma$. See ([P81b]) for 
concrete such examples. For instance, if $\Gamma = \Bbb Z/2\Bbb Z \wr \Bbb Z$ is the lamp-lighter group, then $L(\Gamma)=R$ 
by (Lemma 5.2.3 in [MvN43]; cf. also [C76]) and $H=\Bbb Z$ gives rise to a coarse (abelian) von Neumann subalgebra $L(\Bbb Z)=A \subset R$. 

$2^\circ$ If $B$ is a diffuse tracial von Neumann algebra and $B_0$ is an arbitrary non-trivial tracial 
von Neumann algebra, then $M=B*B_0$ is a II$_1$ factor and $B=B*1\subset M$ is a coarse von Neumann subalgebra, 
by the very definition of free product.

$3^\circ$ If $\Gamma$ is an infinite group, $N_0$ is a non-trivial tracial von Neumann algebra and $\Gamma \curvearrowright N=N_0^{\overline{\otimes} \Gamma}$ 
is the Bernoulli $\Gamma$-action with base $N_0$, then by an argument in ([J81]) it follows that $L(\Gamma)\subset M=N\rtimes\Gamma$ is coarse with infinite multiplicity. 

\proclaim{2.5.4. Proposition} Let $B$ be a diffuse von Neumann subalgebra of the tracial von Neumann algebra $M$. 

$1^\circ$ If $B\subset M$ coarse, then $B'\cap M = \Cal Z(B) \supset  \Cal Z(M)$. 

$2^\circ$ Assume $B\subset M$ is coarse. If  $z_0\in \Cal Z(M)$ is the maximal central projection such that $Bz_0=Mz_0$, 
then $B=Mz_0 + B(1-z_0)$ with $B_1=B(1-z_0) \subset M(1-z_0)=M_1$ being coarse. Moreover,  the representation of $B\overline{\otimes} B^{op}$ 
on $L^2(M\ominus B)=L^2(M_1\ominus B_1)$ factors to a faithful representation of $B_1 \overline{\otimes} B_1^{op}$. 
\endproclaim
\noindent
{\it Proof}. $1^\circ$ If $x\in (B'\cap M) \ominus \Cal Z(B)$, then $L^2(Bx)\subset L^2(M\ominus B)$. 
Thus, $L^2(Bx)$ is a multiple of the coarse bimodule, while at the same time it is a multiple of the trivial $B$-bimodule. 
This forces $x=0$. 

$2^\circ$ Notice first that there exist no projections $q_0, q_1 \in \Cal Z(B)$ such that $L^2(q_0Mq_1)$ is non-zero and finite index  
Hilbert $Bq_0-Bq_1$ bimodule (see also 2.6.3 below). Thus, we may assume $B\subset M$ is so 
that there are no projections $q_0, q_1\in B(1-z_0)$ such that $q_0Mq_1\neq 0$ and $_{Bq_0} L^2(q_0Mq_1)_{Bq_1}$ has finite index. So in order to prove 
the statement we may assume $B\subset M$ itself has the property that $_{Bq_0} L^2(q_0Mq_1)_{Bq_1}$ has infinite index, for any $q_0, q_1\in \Cal P(B)$ 
with $q_0Mq_1\neq 0$. In particular, $M\not\prec_M B$. 

Assume there exists a non-trivial central projection $z\in \Cal Z(B \overline{\otimes} B^{op})$ on which the corresponding representation 
of $B\overline{\otimes} B^{op}$ on $L^2(M\ominus B)$ vanishes. Since $\Cal Z(B \overline{\otimes} B^{op})=\Cal Z(B)\overline{\otimes}\Cal Z(B)^{op}$, 
we can approximate $z$ arbitrarily well in the Hilbert-norm implemented by $\tilde{\tau}=\tau\otimes \tau$ 
by a projection of the form $\Sigma_i z_i J(p_i)J$, where $\{p_i\}_i$ is a finite partition of $1$ with 
projections in $\Cal Z(B)$ and $z_i\in \Cal P(\Cal Z(B))$, $\forall i$. This would imply $\|\Sigma_i p_i x z_i\|_2$ small uniformly 
in $x\in (M\ominus B)_1$. But $M\not\prec_M B$, so there exists a unitary $u\in M$ that's almost orthogonal to $B$. Taking $x=u-E_B(u)$ 
implies $\Sigma_i \tau(p_i)\tau(z_i)$  close to $0$, i.e., $\tilde{\tau}(z)$ arbitrarily close to $0$, a contradiction.  
\hfill $\square$

\vskip.05in
\noindent
{\bf 2.6. Strong malnormality and mixing}. We relate here coarseness with mixing and malnormality properties of subalgebras. 

\vskip.05in
\noindent
{\it 2.6.1. Definition}. A diffuse von Neumann subalgebra $B$ of a  tracial von Neumann algebra $M$ is {\it strongly malnormal} in $M$ 
if any $x\in M$ for which there exists a diffuse abelian von Neumann subgebra $A_0\subset B$ such that 
the Hilbert right $B$-module $L^2(A_0 x B)$ is finite dimensional over $B$, must be contained in $B$. 
With the notation from 2.2 above (cf. Section 1.3 in [P16]), this amounts to the intertwining 
space $\Cal I_M(A_0, B)$ being contained in $B$, for any $A_0\subset B$ diffuse.  
Note that  if $B$ is strongly malnormal and $x\in M$ is so that $\text{\rm dim}(_BL^2(BxA_0))< \infty$, then $x\in B$. 

If $B, B_0\subset M$ are von Neumann subalgebras, then we denote by $w\Cal I_M(B, B_0)$ the 
space of {\it weak intertwiners} from $B$ to $B_0$, i.e.,  of $x\in M$ for which  
there exists a diffuse subalgebra $A_0\subset B$ such $x\in \Cal I_M(A_0, B_0)$. Thus, $B$ strongly 
malnormal in $M$ means that $w\Cal I_M(B, B) \subset B$, or equivalently $w\Cal I_M(B, B)\cap B^\perp=0$. 

Note that  if $B$ is strongly malnormal in $M$ then it is {\it malnormal} in $M$, in the sense of ([IPeP05], [GP14]), i.e., 
if $u\in \Cal U(M)$ satisfies $uBu^*\cap B$ diffuse, then $u\in B$.  

Like in (Definition 1.2.2 in [IPeP05]), given any von Neumann subalgebra $B\subset M$ one can construct by (transfinite) induction the smallest 
von Neumann subalgebra $\tilde{B} \subset M$ that contains $B$ and is strongly malnormal, by considering the strictly increasing family 
of von Neumann algebras $B=B_0 \subset B_1 \subset .... \subset B_\iota=:\tilde{B}$, indexed by the first $\iota$-ordinals such that: $(a)$ for each 
$j < \iota$, one has $w\Cal I(B_j,B_j)\neq B_j$ and $ B_{j+1} $ is the von Neumann algebra generated by $w\Cal I(B_j,B_j)$; $(b)$ if $j\leq \iota$ 
has no predecessor then $B_j=\overline{\cup_{i<j} B_i}$; $(c)$ $w\Cal I(B_\iota, B_\iota)\subset B_\iota$. We call $\tilde{B}$ the 
{\it strongly malnormal cover} (or {\it envelop}) of $B$ in $M$. 
\vskip.05in
\noindent
{\it 2.6.2. Definition}. Let $B\subset M$ be a diffuse von Neumann subalgebra. 
Following (2.9 in [P05b]), we say that $\Cal U(B) \curvearrowright^{\text{\rm Ad}} M$ is {\it mixing relative 
to} $B$ if $\underset{u \rightarrow 0}\to{\lim}  \|E_B(xuy)\|_2 = 0,$ 
for all $x, y\in M\ominus B$, where the limit is over $u\in \Cal U(B)$ tending weakly to $0$. We then also say that 
$B\subset M$ is {\it mixing}. In this same vein, if $B_0\subset M$ is another von Neumann subalgebra and 
we have $\lim_u \|E_{B_0}(xuy)\|_2=0$, for any $x, y\in M$, 
where the limit is over unitaries $u\in B$ that tend to $0$ in the weak operator topology, then we say that $B, B_0\subset M$ is a 
{\it mixing pair} of subalgebras.

\proclaim{2.6.3. Proposition} $(a)$ If $B \subset M$ is a diffuse von Neumann subalgebra of the tracial von Neumann $M$, 
then $B$ coarse $\Rightarrow$ $B$ mixing $\Rightarrow$ $B$ strongly malnormal.  

$(b)$ If $B\subset N$ and $N\subset M$ are coarse $($respectively mixing, resp. strongly malnormal$)$, then $B\subset M$ 
is coarse $($resp. mixing, resp. strongly malnormal$)$. 

$(c)$  Let  $B, B_0 \subset M$ be diffuse von Neumann subalgebras of the tracial von Neumann algebra $M$. 
If $B, B_0$ is a coarse pair then $B, B_0$ is a mixing pair. Also, if $B, B_0$ is a mixing pair, then  $w\Cal I_M(B, B_0)=0$.
\endproclaim
\noindent
{\it Proof}. $(a)$ By using the polarization trick, showing  that $B\subset M$ is mixing is equivalent to showing that 
 $\underset{u \rightarrow 0}\to{\lim}  \|E_B(xux^*)\|_2 = 0,$ 
for all $x\in M\ominus B$ with $\|x\|_2=1$, a condition that's equivalent to $\underset{u \rightarrow 0}\to{\lim}(\sup \{|\tau(bxux^*)| \mid b\in (B)_1\})=0$. 

If one denotes by $\tilde{\varphi}$ the state  on $B\vee B^{op}$ implemented by $\hat{x}\in L^2M\ominus L^2B$, 
then $B$ coarse in $M$ implies that $\tilde{\varphi}$ implements a normal state on $B\overline{\otimes} B^{op}$.  
So $\tilde{\varphi}$ is of the form $\tilde{\tau}( \cdot \ \tilde{b})$, for some $\tilde{b}\in L^1(B\overline{\otimes} B^{op}, \tilde{\tau})_+$ 
with $\tilde{\tau}(\tilde{b})=1$, where $\tilde{\tau}=\tau\otimes \tau$. But for any such 
normal state $\tilde{\varphi}$ on $B\overline{\otimes} B^{op}$ one has $\underset{\tilde{y} \rightarrow 0}\to{\lim} \tilde{\varphi}(\tilde{y})=0$, where the limit is 
taken over $\tilde{y}\in (B\overline{\otimes} B^{op})_1$ tending weakly to $0$. Since $u\in \Cal U(B)$ tending weakly to $0$ implies 
$b \otimes u^{op}$ tends weakly to $0$ uniformly in $b\in (B)_1$, it follows that 
$\underset{u \rightarrow 0}\to{\lim}(\sup \{|\tilde{\varphi}(b\otimes u^{op})| \mid b \in (B)_1\})=0$, 
where the limit is over $u\in \Cal U(B)$ tending weakly to $0$.  By taking into account that $\tilde{\varphi}(b \otimes u^{op})=\tau(bxux^*)$, 
this shows that $B\subset M$ coarse implies $B\subset M$ mixing. 

The proof that $B\subset M$ mixing implies $B\subset M$ strongly malnormal is exactly as the proof of (3.1 in [P03]), so we leave the details to the reader. 

\vskip.05in

Part $(b)$ is trivial, by the definitions, while the proof of  part $(c)$ is very similar to the proof of part $(a)$ above, so we leave the details to the reader. 

\hfill 
$\square$

\vskip.05in
\noindent
{\it 2.6.4. Examples} $1^\circ$ The examples of coarse subalgebras and course pairs in Example 2.5.3.1$^\circ$, 
arising from subgroups  $H, H_0 \subset \Gamma$, have been much exploited in [P81b], where one implicitly gives a proof of the above proposition and 
apply this  to provide many  examples of MASAs in group factors, with calculable normalizers.  

$2^\circ$ The coarse inclusions $B=B*1\subset B * B_0=M$ in Example 2.5.3.2$^\circ$ follow mixing and strongly malnormal by Proposition 2.6.3 above. 
Note that both these properties are implicitly proved in ([P90], [IPeP05]).  

$3^\circ$ The example $2.5.3.3^\circ$ of coarse inclusion $L(\Gamma) \subset M = N\rtimes \Gamma$, where $\Gamma \curvearrowright N=N_0^{\overline{\otimes}\Gamma}$ 
is the Bernoulli $\Gamma$-action with base $N_0\neq \Bbb C$, follows 
mixing and strongly malnormal by Proposition 2.6.3 above. One can in fact show quite easily 
that if $\Gamma \curvearrowright (N,\tau)$ is a a free mixing action, then $L(\Gamma) \subset M=N\rtimes \Gamma$  
is mixing. Note that if 
$\Gamma=\Bbb Z$ and $N_0$ is abelian, then $M$ is isomorphic to the hyperfinite II$_1$ factor $R$. Thus, it follows that 
$A=L(\Bbb Z)\subset R$ is a coarse MASA, which is also strongly malnormal and mixing. Moreover, it is easy to see that the inclusion 
$A \vee A^{op} \subset  (A \vee A^{op})'\cap \Cal B(L^2R \ominus L^2A)$ is isomorphic to 
$A \overline{\otimes} A^{op} \subset (A \overline{\otimes} A^{op}) \otimes \Cal B(\ell^2 \Bbb N)$, implying that the MASA $A\subset R$ has infinite multiplicity (Pukaszky invariant equal to $\infty$).

\vskip.05in 

\noindent
{\bf 2.7. A criterion for $R$-bimodules to be coarse}. We end this section with a criterion for commuting normal representations $R, R^{op}$  
of the hyperfinite II$_1$ factor on the same Hilbert space $\Cal H$ to generate the II$_1$ tensor product 
$R\overline{\otimes} R^{op}$ in $\Cal B(\Cal H)$, i.e., for the Hilbert 
$R$-bimodule $\Cal H$ to be a multiple of the coarse $R$-bimodule. This will be an immediate  consequence of the following:   

\proclaim{2.7.1. Lemma}  Let $Q$ be a tracial von Neumann algebra represented normally and faithfully 
on a separable Hilbert space $\Cal H$. Let $\Cal B\subset Q'\cap \Cal B(\Cal H)$ be a  UHF algebra, 
obtained as the $C^*$-inductive limit of matrix factors $\Cal B_n\simeq \Bbb M_{k_n}(\Bbb C)$ 
with $k_n | k_{n+1}$, $\forall n$, and $\tau$ its unique trace state. Denote by $\Cal Q$ the $C^*$-inductive limit of $\Cal Q_n=Q \vee \Cal B_n$. 
Let  $\{\xi_n\}_n\subset \Cal H$ be 
a sequence of unit vectors in $\Cal H$ that's dense in the set of unit vectors of $\Cal H$. The following conditions are equivalent 

\vskip.05in
$(a)$ The von Neumann algebra $\Cal Q''=\overline{\cup_n \Cal Q_n}^{wo}$ is the tracial von Neumann 
algebra $Q\overline{\otimes} R$; 

\vskip.03in
$(b)$ Given any unit vector $\xi \in \Cal H$, the vector state $\varphi$ implemented by $\xi$ satisfies 
$\underset{n\rightarrow \infty}\to{\lim} \|\tau_{|\Cal B_n'\cap \Cal B} \otimes \varphi_{|Q} - \varphi_{|\Cal B_n'\cap \Cal Q}\|=0$. 

\vskip.03in
$(c)$ The vector state $\varphi_m$ implemented by $\xi_m$ satisfies 
$\underset{n\rightarrow \infty}\to{\lim} \|\tau_{|\Cal B_n'\cap \Cal B}\otimes {\varphi_m}_{|Q}-{\varphi_m}_{|\Cal B_n'\cap \Cal Q}\|=0$, $\forall m$.

\endproclaim
\noindent
{\it Proof}. This is a straightforward  consequence of  criteria in [Po67]. We leave the details as an exercise. 
\hfill 
$\square$

\proclaim{2.7.2. Corollary} Let $R, R^{op} \subset \Cal H$ be commuting normal representations of the hyperfinite $\text{\rm II}_1$ factor and its 
opposite on the separable Hilbert space $\Cal H$. Let $B_n\subset R$ be an increasing sequence of matrix factors such that $(\cup B_n)''=R$ 
and $\{\xi_n\}_n \subset \Cal H$ a dense subsequence of the unit sphere of $\Cal H$. Let $\Cal R$ denote 
the $^*$-algebra $\Cal R=\text{\rm Alg}(R, R^{op})$ generated by $R, R^{op}$ in $\Cal B(\Cal H)$ which we identify 
in the usual way with $R\otimes R^{op}$. Denote by $\tilde{\tau}$  the trace state on $\Cal R$ 
defined by $\tilde{\tau}(x \otimes y^{op})=\tau(x)\tau(y)$, $\forall x, y \in R$, and by $\varphi_m$ the vector state implemented by $\xi_m$, $m\geq 1$.  
The von Neumann algebra $\Cal R''\subset \Cal B(\Cal H)$ is a $\text{\rm II}_1$ factor  $(\simeq R\overline{\otimes} R^{op})$ if and only if  
$$
\lim_{n \rightarrow \infty} \|\tilde{\tau}_{|(B_n \otimes B_n^{op})'\cap \Cal R}-{\varphi_m}_{|(B_n \otimes B_n^{op})'\cap \Cal R}\|=0, \forall m. 
$$
\endproclaim
\noindent
{\it Proof}. Immediate by Lemma 2.7.1 applied to $Q=\Bbb C$ and $\Cal B_n=B_n \vee B_n^{op}$.   
\hfill 
$\square$ 

\heading 3. A technical lemma  \endheading 

In this section we prove a key technical result needed in the proof of Theorem 1.1.  The proof uses 
the incremental patching technique, in a manner similar to ([P92], [P13a], [P13b], [P17]).  

\proclaim{3.1. Lemma} Let $M$ be a $\text{\rm II}_1$ factor, 
$Q\subset M$ a von Neumann subalgebra,  $P \subset M$ an irreducible subfactor such that $P\not\prec_M Q$ 
and $P_0\subset P$ a finite dimensional subfactor. 
Given any finite sets $F=F^*\subset (M\ominus P_0)_1$, $1\in F'={F'}^*\subset (M)_1$, and any $\delta_0 >0$, there 
exists a unitary element $v_0\in P_0'\cap P$ such that 

$$
\|E_Q(x_0v_0xv_0^*x_0^*)\|_2^2 \leq \delta_0, \forall x\in F, x_0\in F'; \tag 3.1.1
$$
$$ 
|\tau(v_0x_1v_0^*x_2v_0x_3v_0^*x_4)| \leq \delta_0, \forall x_i\in F. \tag 3.1.2
$$ 
$$ 
|\tau(v_0xv_0^*y)| \leq \delta_0, \forall x, y \in F.  \tag 3.1.3
$$ 

\endproclaim

\noindent
{\it Proof}.  Let $\omega$ be a non-principal ultrafilter on $\Bbb N$ 
and denote by $\Cal M=\langle M^\omega, e_{Q^\omega \rangle}$ the semifinite von Neumann algebra 
associated with the basic construction for $Q^\omega \subset M^\omega$. 
Thus, $\Cal M = \overline{\text{\rm sp}M^\omega e M^\omega}^w \subset \Cal B(L^2M^\omega)$, 
where $e=e_{Q^\omega} \in \Cal B(L^2M^\omega)$. 

Fix $\delta  >0$ such that $\delta < \delta_0$. Denote by ${\mycal W}$ the set of partial isometries $v\in P_0'\cap P^\omega =(P_0'\cap P)^\omega$  
with the property that $vv^*=v^*v$ and which 
satisfy the conditions: 
$$
 \|E_{Q^\omega}(x_0vxv^*x_0^*)\|_2^2 \leq \delta \tau(v^*v),  \forall x\in F, x_0\in F'  \tag a
$$
$$
|\tau(vx_1v^*x_2vx_3v^*x_4)|\leq \delta \tau(v^*v), \forall x_i \in F \tag b 
$$
  $$
|\tau(v^*xvy)| \leq \delta\tau(v^*v), \forall x, y\in F  \tag c
 $$ 
$$
 E_{P_0}(vv^*F) = 0, \tag d
 $$
 $$
 E_{P_0}(FvF)=0, \tag e
 $$ 
  $$
 E_{P_0}(vv^*FvF)=0. \tag f
 $$

We endow ${\mycal W}$ with the order $\leq$  
in which $v_1 \leq v_2$ if $v_1=v_2v_1^*v_1$. 
$({\mycal W}, \leq)$ is then clearly inductively ordered and we let $v\in {\mycal W}$ be a maximal element. 

Assume $\tau(v^*v) < 1$ and denote $p = 1 - v^*v$.  
Notice that $p\in P_0'\cap P^\omega$ and that by $(c), (d), (e)$ we have 
$E_{P_0}(pF)=0$ and $E_{P_0}(pFvF)=0$. Since the uniqueness of trace preserving expectation onto $P_0$ implies 
that for a unitary element $u\in P_0'\cap M^\omega$ and $y\in M^\omega$ we have 
$E_{P_0}(y)=u^*E_{P_0}(uyu^*)u= E_{P_0}(uyu^*)$, it follows that for any $x\in F$ and $u\in \Cal U(P_0'\cap P^\omega)$ 
we have $E_{P_0}(upx)=E_{P_0}(u(pxu)u^*)=E_{P_0}(pxu)$. By writing $p$ as a linear combinations between $u=1$ and $u=2p-1$ 
this implies that $E_{P_0}(pF)=E_{P_0}(pFp)$ and thus $E_{P_0}(pFp)=0$ as well. Similarly, $E_{P_0}(pFvFp)=0$. 

Let $w$ be a
partial isometry in $p(P_0'\cap P^\omega)p$ with $w^*w=ww^*$ and denote $u=v+w$. Then $u$ is a partial isometry in $P_0'\cap P^\omega$ with $u^*u=uu^* 
\in P_0'\cap P^\omega$. 
We will show that one can make an appropriate choice 
$w\neq 0$ such that $u = v + w$ lies in ${\mycal W}$. This will contradict the maximality of $v$, thus showing that $v$ must be a unitary element. 
We will construct  the partial isometry $w$ by first choosing its support $q=ww^*=w^*w$, then choosing the ``phase $w$'' above $q$.

In order to get estimates on $(a)$, note that by writing  
$eux^*u^*x_0^*ex_0uxu^*$ as $e(v+w)x^*(v+w)^*x_0^*ex_0(v+w)x(v+w)^*$ and developing into the sum of 16 terms, we get 
$$
\|E_{Q^\omega}(x_0uxu^*x_0^*)\|_2^2=Tr(ex_0ux^*u^*x_0^*ex_0uxu^*x_0^*)  \tag 1a
$$
$$
\leq Tr(ex_0vx^*v^*x_0^*ex_0vxv^*x_0^*) + \Sigma_{1,a} +\Sigma_{2,a}+\Sigma_{3,a}+ \Sigma_{4,a}, 
$$
where $\Sigma_{i,\text{\rm a}}$ denotes the sum of  the absolute value of terms having $i$ appearances of elements from $\{w, w^*\}$, $1\leq i \leq 4$. 
Thus, there are four terms in $\Sigma_{1,a}$, six in $\Sigma_2(a)$, four in $\Sigma_3(a)$, and one in $\Sigma_4(a)$. 

Similarly, in order to estimate $(b)$, by developing $\tau(ux_1u^*x_2ux_3u^*x_4)=\tau((v+w)x_1(v+w)^*x_2(v+w)x_3(v+w)^*x_4)$  into a sum of 16 terms  we get 
$$
|\tau(ux_1u^*x_2ux_3u^*x_4)| \leq |\tau(vx_1v^*x_2vx_3v^*x_4)| +  \sum_{i=1}^4 \Sigma_{i,\text{\rm b}}   \tag 1b
$$
where $\Sigma_{i,\text{\rm b}}$ denotes the sum of absolute value of the 
terms $\tau(y)$ with $y$ having $i$ appearances of elements from $\{w, w^*\}$, $1\leq i \leq 4$. Also, for $(c)$ we have 
$$
 |\tau(u^*xuy)| \leq |\tau(v^*xvy)| + |\tau(wxv^*y)| + |\tau(vxw^*y)|+ |\tau(wxw^*y)|, \forall x, y \in F  \tag 1c
$$

At the same time, in order for $(d), (e), (f)$ to be satisfied, we need to have: 

$$
E_{P_0}(ww^*F)=0 \tag 1d 
$$
$$
E_{P_0}(FwF)=0 \tag 1e
$$
$$
E_{P_0}(ww^*FvF)=E_{P_0}(ww^*FwF)=E_{P_0}(vv^*FwF)=0 \tag 1f 
$$

Let us first  estimate the terms $|Tr(X)|$ in $(1a)$  
with $X$ containing a pattern of the form $... ex_0wxw^*x_0^*e...$, or $...ex_0wx^*w^*x_0^*e...$, for given $x\in F$, $x_0\in F'$. There are 
seven such terms : the one in $\Sigma_{4,\text{\rm a}}$, all four in $\Sigma_{3,\text{\rm a}}$ 
and two in $\Sigma_{2,\text{\rm a}}$.  We denote by $\Sigma'_{\text{\rm a}}$ the sum 
of these terms. Note that for each such $X$ we have $|Tr(X)|=|Tr(wxw^*y'ey)|$ for some $y,y' \in (M^\omega)_1$. 
Thus, by applying the Cauchy-Schwartz inequality and taking into account the definition of $Tr$, we get the estimate 
$$
|Tr(X)|=|Tr(wxw^*y'ey)| \tag 2a
$$
$$
\leq (Tr(e{y'}^*wx^*w^*wxw^*y'e))^{1/2}(Tr(q'ey^*yeq'))^{1/2} \leq \|qxq\|_2\|q\|_2, 
$$
where $q'$ is the left support of $y'w$, which thus has trace $\leq \tau(q)$, implying  that $Tr(q'ey^*yeq')\leq Tr(qeq)=\tau(q)$. We have also used   
that $Tr(e{y'}^*wx^*w^*wxw^*y'e)=Tr(e{y'}^*wx^*qxw^*y'e) \leq \tau(wx^*qxw)=\tau(qx^*qxq)$. 

Similarly, the seven terms $|\tau(y)|$  in (1b) with $y$ containing a pattern of the form 
$... wx_jw^*...$, or $...w^*x_iw...$ (namely, the one in $\Sigma_{4,\text{\rm b}}$, all four in $\Sigma_{3,\text{\rm b}}$ 
and two in $\Sigma_{2,\text{\rm b}}$) are  majorized by 
$$
|\tau(y)| \leq \|qx_jq\|_2 \|q\|_2. \tag 2b
$$
In addition, since $pFvFp \perp P_0p$, for the remaining two terms  $y=wx_1v^*x_2wx_3v^*x_4$, $y=vx_1w^*x_2vx_3w^*x_4$ in $\Sigma_{2,\text{\rm b}}$,  we have 
$$
|\tau(y)| \leq \|q(x_1vx_2)q\|_2 \|q\|_2.  \tag 2b'
$$

As for $(1c)$, for the only term containing both $w, w^*$ we have the estimate
$$
 |\tau(wxw^*y)|\leq \|qxq\|_2 \|qyq\|_2 \tag 2c
 $$  

By  (Lemma 2.3; cf. 1.4 in [P92], or 4.3 in [P13b]), since $pFp, pFvFp$ are perpendicular to $P_0p$, the subfactor  $p(P_0'\cap P^\omega)p$ of  the II$_1$ factor $pM^\omega p$ 
contains a diffuse abelian subalgebra that's $2$-independent to $pFp, pFvFp$ relative to $P_0p$   
with respect to the trace state $\tau( \cdot )/\tau(p)$  on $pM^\omega p$. This implies that  
there exists a projection $q\in p(P_0'\cap P^\omega)p$ of trace $\tau(q)=\delta^2 \tau(p)^2/12^2$ such that $E_{P_0p}(q(pFp))=0$, $E_{P_0p}(q(pFvFp))=0$ 
and $\|qzq\|_2^2/\tau(p)=(\tau(q)/\tau(p))^2\tau(z^*z)/\tau(p)$, for all $z\in pFp\cup pFvFp$. 
 
Since $q\leq p$, it follows that for each  $x\in F$ one has 
$$\|qxq\|^2_2 \leq (\delta^4\tau(p)^2/12^4) \tau(x^*x)\leq\delta^2 \tau(q)/12^2.$$ Thus,    
$\|qxq\|_2\leq \delta \tau(q)^{1/2}/12$, $\forall x\in F$. Hence,  
for this choice of $q$, the  right hand side term in (2a) will be majorized by $\delta \tau(q)/12$.  
By summing up over the seven terms in $\Sigma'_{\text{\rm a}}$,  we get 
$$
\Sigma'_{\text{\rm a}} \leq 7\delta \tau(q)/12 \tag 3a
$$

Similarly for the eleven terms in $\Sigma_{2,\text{\rm b}}, \Sigma_{3,\text{\rm b}}, \Sigma_{4,\text{\rm b}}$ we get:
$$
\Sigma_{2,\text{\rm b}} +\Sigma_{3,\text{\rm b}}+\Sigma_{4,\text{\rm b}} \leq 11\delta \tau(q)/12. \tag 3b
$$ 
while for the single term in $(2c)$ we get
$$
 |\tau(wxw^*y)|\leq \|qxq\|_2 \|qyq\|_2 \leq \delta^2 \tau(q)/144   \tag 3c
 $$  

We will now estimate the sum $\Sigma''_{\text{\rm a}}$ of  the terms $|Tr(X)|$ with $X$ running over the remaining four terms in $\Sigma_{2,\text{\rm a}}$, 
the sum $\Sigma_{1,\text{\rm a}}$ of the four terms $|Tr(X)|$ with $X$ having only one occurrence of $w, w^*$,  
the sum $\Sigma_{1,\text{\rm b}}$ of the four terms $|\tau(y)|$ in (1b) with $y$ having only one occurrence of $w, w^*$, and the 
two terms in $(1c)$ with just one occurence of $w, w^*$, while at the 
same time taking care of the conditions $(1c)-(1f)$ (of which $(1c)$ and the first equality in  $(1f)$ are already satisfied by 
the choice of $q=ww^*$). We will do this by making an appropriate choice 
of the ``phase $w$'' above the support projection $q$, which is  fixed. 

Note that all elements entering in the sums $\Sigma_{1,\text{\rm a}}$, $\Sigma_{1,\text{\rm b}}$ 
and the terms $|\tau(wxv^*y)|$, $|\tau(vxw^*y)|$ in $(1c)$ are of the form $|\tau(wz)|$, where $z$ belongs to a finite 
set $E\subset (qM^\omega q)_1$. 
Let $\{e_{kl}\}_{k,l}$ be matrix units for $P_0$ and let $\Cal F\subset (qM^\omega q)_1$ denote the finite set $q(\cup_{k,l} Fe_{kl}pF \cup E)q$. 
By results in ([P92], [P16]), there exists a hyperfinite subfactor $R\subset q(P_0'\cap P^\omega)q$ such that 
$E_{R'\cap qM^\omega q}(z')=E_{P_0q}(z')$, $\forall z'\in \Cal F$.

Since $P_0q$ and $R'\cap qM^\omega q$ are $\tau$-independent, if we denote $\tau_q$ the normalized trace $\tau( \ )/\tau(q)$ on $qM^\omega q$ then 
for each unitary element  $w\in N:=R'\cap q(P_0'\cap P^\omega) q$ and $z\in \Cal F$ we have 
$$
|\tau(wz)|/\tau(q) = |\tau_q(wz)| = |\tau_q(E_{R'\cap qM^\omega q}(wz))| = |\tau_q(wE_{R'\cap qM^\omega q}(qzq))| 
$$
$$
= |\tau_q(w)| |\tau_q(E_{R'\cap qM^\omega q}(qzq))|= |\tau(w)| |\tau(E_{P_0q}(qzq))|/\tau(q). 
$$
Since $\|z\|\leq 1$, this implies that for any $w\in \Cal U(N)$ we have 

$$
|\tau(wz)|\leq |\tau(w)|, \forall z\in \Cal F, \tag 4
$$
$$
\Sigma_{1,\text{\rm a}} \leq 4 |\tau(w)|, \tag 4a
$$
$$
\Sigma_{1,\text{\rm b}} \leq 4 |\tau(w)|  \tag 4b
$$ 
$$
|\tau(wxv^*y)| +  |\tau(vxw^*y)| \leq 2|\tau(w)| \tag 4c
$$

At this point,  it is convenient to enumerate the elements in $F=\{x_1, ..., x_n\}$, $F'=\{y_1, ..., y_n\}$ 
(we may clearly assume $|F|=|F'|$). For each $1\leq i, k \leq n$ we have  
$$
\Sigma''_{\text{\rm a}}=|Tr(y_k^*ey_kwx_i^*v^*y_k^*ey_kvx_iw^*)|+|Tr(y_k^*ey_kvx_i^*w^*y_k^*ey_kwx_iv^*)| \tag 5a 
$$
$$
+|Tr(y_k^*ey_kwx_i^*v^*y_k^*ey_kwx_iv^*)|+|Tr(y_k^*ey_kvx_i^*w^*y_k^*ey_kvx_iw^*)|
$$
$$
=|Tr(w^*ewY_{1,i,k}|+|Tr(w^*ewY_{2,i,k})|
$$
$$
+|Tr(wY_{3,i,k}wY_{4,i,k})|+|Tr(w^*Y_{5,i,k}w^*Y_{6,i,k})|
$$
where each one of the terms $Y_{j,i,k}$ depends on $x_i\in F, y_k \in F'$ and 
belongs to the set  $S_0$ $:=q((M^\omega)_1e(M^\omega)_1)q$ $\subset qL^2(\Cal M, Tr)q$. 

Note that, as $1\leq i, k \leq n$,  the number 
of possible indices $(j,i,k)$ in  $(5a)$ is $4n^2$. There are $2n^2$ terms of the form $|Tr(w^*ewY)|$, 
$n^2$ terms of the form $|Tr(wXwY)|$ and $n^2$ terms of the form $|Tr(w^*Xw^*Y)|$, which by 
using the fact that $|Tr(w^*Xw^*Y)|=|Tr(wX^*wY^*)|$ we can view as $n^2$ additional terms of the form 
$|Tr(wXwY)|$. In all this, the elements $X, Y$ belong to $S_0\subset qL^2(\Cal M, Tr)q$, and are thus 
bounded in operator norm by $1$ and are supported (from left and right) by projections of trace $Tr$ majorized by $1$.

Recall that we are under  the assumption 
$P\not\prec_M Q$. By Lemma 2.1, this implies $R'\cap (P_0'\cap qP^\omega q) \not\prec_{M^\omega} Q^\omega$.  
Thus,  $N'\cap q\Cal M q$ contains no finite non-zero projections of $\Cal M=\langle M^\omega, Q^\omega \rangle$. 

To estimate the terms in $\Sigma''_{\text{\rm a}}$ (and at the same time $\Sigma_{1,\text{\rm a}}$, $\Sigma_{1,\text{\rm b}}$), we will prove the following:

\vskip.05in 
\noindent
{\it Fact} 1. For any $\alpha >0$ and any two $m$-tuples of elements $(Z_1, ..., Z_m)$, 
$(Z_1', ..., Z_m')$ in $S_0 \cap \Cal M_+$, 
there exists a unitary element $w\in N$ such that $|\tau(w)|\leq \alpha/4$ and 
$$
\sum_{i=1}^m Tr(w^*Z_iwZ_i')\leq \alpha. 
$$

\vskip.05in

To prove this, let $\Cal H$ denote the Hilbert space $L^2(q\Cal Mq, Tr)^{\oplus m}$ and note that we have a unitary 
representation $\Cal U(N)\ni w \mapsto \pi(w) \in \Cal U(\Cal H)$, which on an $m$-tuple 
$X=(X_i)_{i=1}^{m}\in \Cal H$ acts by $\pi(w)(X)=(w^*X_iw)_i$.  

Now note that this representation has no  
(non-zero) fixed point. Indeed, for if $X\in \Cal H$ satisfies $\pi(w)(X)=X$, $\forall w\in \Cal U(N)$, 
then on each component $X_i\in L^2(q\Cal Mq, Tr)$ of $X$ we would have $w^*X_iw=X_i$, $\forall w$. 
Thus $X_iw=wX_i$ and since the unitaries of $N$ span linearly the algebra $N$, 
this would imply $X_i\in N'\cap L^2(q\Cal Mq, Tr)$. Hence, $X_i^*X_i\in N'\cap L^1(q\Cal Mq, Tr)$ and therefore all spectral 
projections of $X_i^*X_i$ corresponding to intervals $[t, \infty)$ with $t>0$ would be projections 
of finite trace in $N'\cap q\Cal Mq$, forcing them all to be equal to $0$. Thus, $X_i=0$ for all $i$. 

With this in mind, denote by $K_Z\subset \Cal H$ the weak closure of the convex hull of the set $\{\pi(w)(Z) \mid w\in \Cal U(N)\}$, 
where $Z=(Z_1, ..., Z_m)$ is viewed as an element in $\Cal H$. Since $K_Z$ is bounded and weakly closed, it is weakly compact, 
so it has a unique element $Z^0\in K_Z$ of minimal norm $\| \ \|_{2,Tr}$. Since $K_Z$ is invariant to $\pi(w)$ and $\|\pi(w)(Z^0)\|_{2,Tr}=\|Z^0\|_{2,Tr}$, 
it follows that $\pi(w)(Z^0)=Z^0$. But we have shown that $\pi$ has no non-zero fixed points, and so $0=Z^0\in K_Z$. 

Let us deduce from this that if $Z=(Z_i)_i, Z'=(Z'_i)_i$ are the two $m$-tuples of positive elements in $S_0$, then 
we can find $w\in \Cal U(N)$ such that the inequality in {\it Fact} 1 holds true. Indeed, for if there would exist $\alpha>0$ such that $\sum_i Tr(\pi(w)(Z_i)Z_i')\geq \alpha$, 
$\forall w\in \Cal U(N)$, then by taking convex combinations and weak closure, one would get $0=\langle Z^0, Z'\rangle \geq \alpha$,  
a contradiction. 

Note that by taking one of the $i$  elements $Y_i, Y_i'$ to be equal to $e$, 
one can get $w\in \Cal U(N)$ to also satisfy  $|\tau(w)|^2\leq \alpha^2/16$.  This finishes the proof of  {\it Fact} 1. 

We will now use {\it Fact} 1 to prove: 

\vskip.05in 
\noindent
{\it Fact} 2.  Given any $m$-tuples $(X_i)_i, (Y_i)_i, (X'_i)_i,$ $(Y'_i)_i \in S_0^m$ (not necessarily 
having positive operators as entries) and any $\alpha>0$, there exists $w\in \Cal U(N)$ such that 
$|\tau(w)| \leq \alpha/4$, $\sum_{i=1}^m |Tr(w^*X_iwX'_i)|\leq \alpha$, 
$\sum_{i=1}^m |Tr(wY_iwY_i')|\leq \alpha$. 
\vskip.05in

Indeed, because if we denote by $e_i$ the left support of $X'_i$ and 
$f_i$ the left support of $Y'_i$, then by the Cauchy-Schwartz inequality we simultaneously have for all $i$ the estimates 
$$
|Tr(w^*X_iwX'_i)|^2\leq Tr(w^*X_i^*X_iwX'_i{X'_i}^*)Tr(e_i) 
\leq Tr(w^*X_i^*X_iwX'_i{X'_i}^*), 
$$
and respectively 
$$
|Tr(wY_iwY_i')|^2\leq Tr(w^*Y_i^*Y_iwY'_i{Y'_i}^*)Tr(f_i)
 \leq Tr(w^*Y_i^*Y_iwY'_i{Y'_i}^*). 
$$
Since all $X_i^*X_i, X'_i{X'_i}^*, Y_i^*Y_i, Y'_i{Y'_i}^*$ are positive elements in $S_0$, we can now apply the {\it Fact} 1 
to deduce that  there exist $w\in \Cal U(N)$ such that  $|\tau(w)|\leq \alpha/4$,  
$\sum_{i=1}^m|Tr(w^*X_iwX'_i)|\leq \alpha$, $\sum_{i=1}^m |Tr(wY_iwY_i')|\leq \alpha$. This ends the proof of {\it Fact} 2. 

\vskip.05in 

For each $n\geq 1$, we apply  {\it Fact} 2 to $\alpha=2^{-n-1}$ to get a partial isometry $w_n\in P_0'\cap P^\omega$ 
of support $w_nw_n^*=w_n^*w_n=q$ such that if we denote by 
$\Sigma_{1,a}(w_n)$, $\Sigma_{1,b}(w_n)$, $\Sigma''_{\text{\rm a}}(w_n)$ the values of $\Sigma_{1,a}$, $\Sigma_{1,b}$, $\Sigma''_{\text{\rm a}}$, 
obtained by plugging in $w_n$ for $w$ in the inequalities $(4a)$, $(4b)$, $(5a)$, respectively, then we have: 

$$
|\tau(w_nz)| \leq 2^{-n-1}, \forall z\in \Cal F;
$$
$$
\Sigma_{1,a}(w_n) \leq 2^{-n-1}, \Sigma_{1,b}(w_n) \leq 2^{-n-1}, \Sigma''_{\text{\rm a}}(w_n)\leq 2^{-n-1}. 
$$

Let $q=(q_n)_n$, with $q_n\in \Cal P(P_0'\cap P)$ and $w_n=(w_{n,k})_k$ with $w_{n,k}^*w_{n,k}=w_{n,k}w_{n,k}^*=q_k$, $\forall k$. 
For each $k\geq 1$, we denote $\Cal F_k\subset (q_kMq_k)_1$ the set of all $k$'th entries of elements $z=(z_k)_k \in \Cal F\subset (qM^\omega q)_1$. We also denote by 
$\Sigma_{1,a}(w_n)_k$ (resp. $\Sigma_{1,b}(w_n)_k$, $\Sigma''_{\text{\rm a}}(w_n)_k$) the sum obtained at  the $k$'th level of  
$\Sigma_{1,a}(w_n)$ (resp. $\Sigma_{1,b}(w_n)$, $\Sigma''_{\text{\rm a}}(w_n)$). Thus, we have 
$$
\underset{k\rightarrow \omega}\to\lim |\tau(w_{n,k}z_k)| \leq 2^{-n-1}, \forall z_k\in \Cal F_k;
$$
$$
\underset{k\rightarrow \omega}\to\lim \Sigma_{1,\text{\rm a}}(w_n)_k  \leq 2^{-n-1}; 
$$
$$
\underset{k\rightarrow \omega}\to\lim \Sigma_{1,\text{\rm b}}(w_n)_k  \leq 2^{-n-1}; 
$$
$$
\underset{k\rightarrow \omega}\to\lim\Sigma''_{\text{\rm a}}(w_n)_k  \leq 2^{-n-1}.  
$$

Denote by $V_n$ the set of all $k\in \Bbb N$ such that   $|\tau(w_{n,k}z_k)|<2^{-n}, \forall z_k\in \Cal F_k$, 
$\Sigma_{1,a}(w_n)_k<2^{-n}$, $\Sigma_{1,b}(w_n)_k  \leq 2^{-n}$, $\Sigma''_{\text{\rm a}}(w_n)_k  \leq 2^{-n}$. 
Note that $V_n$ corresponds to an open closed neighborhood of $\omega$ in $\Omega$, 
under the identification $\ell^\infty\Bbb N = C(\Omega)$. Let now $W_n$, $n\geq 0$, be defined recursively as follows: 
$W_0=\Bbb N$ and $W_{n+1}=W_n \cap V_{n+1}\cap \{k \in \Bbb N \mid k > \min W_n\}$. 
Note that, with the same identification as before, $W_n$ is a strictly decreasing sequence of neighborhoods of $\omega$ in $\Omega$. 

Define $w=(w'_m)_m$, where $w'_k=w_{m, k}$ for $k\in W_{m-1}\setminus W_m$. It is then easy to see 
that $w$ is a  partial isometry in $P_0'\cap P^\omega$ with $ww^*=w^*w=q$ and that we have 
$\tau(wz)=0$, $\forall z\in \Cal F$, $\Sigma_{1,\text{\rm a}}(w)=0$, $\Sigma_{1,\text{\rm b}}(w)=0$, $\Sigma''_{\text{\rm a}}(w)=0$. By taking into account the definition of $\Cal F$, 
it is easy to see that the first of these 
conditions implies that $FwF\perp P_0$, $uu^*FwF \perp P_0$, where $u=v+w$. By (1b), (3b) and $\Sigma_{1,\text{\rm b}}(w)=0$, it follows 
that   $|\tau(ux_1u^*x_2ux_3u^*x_4)|\leq \delta \tau(uu^*)$, $\forall x_i \in F$, while from $(2c), (3c)$ we have 
$|\tau(uxu^*y)|\leq \delta\tau(uu^*)$, $\forall x, y\in F$. In turn, from (1a),  (3a) and the fact that $\Sigma_{1,\text{\rm a}}(w)=0$, 
$\Sigma''_{\text{\rm a}}(w)=0$, it follows that $\|E_{Q^\omega}(x_0uxu^*)\|_2^2 \leq \delta \tau(uu^*)$, $\forall x\in F, x_0\in F'$. This shows that $u\in \Cal W$, 
while  $u\geq v$, $u\neq v$, contradicting the maximality of $v$

This shows that $v$ must be a unitary element. Thus, if we represent $v\in P_0'\cap P^\omega$ as a sequence of unitary 
elements $(v_n)_n$ in $P_0'\cap P$, then we have 
$$
\lim_{n\rightarrow \omega} \|E_Q(x_0v_nxv_n^*x_0^*)\|_2^2 = \|E_Q(x_0vxv^*x_0^*)\|_2^2 \leq \delta < \delta_0, 
$$
$$
\lim_{n\rightarrow \omega} |\tau(v_nx_1v_n^*x_2v_nx_3v_n^*x_4)| \leq \delta < \delta_0, 
$$
$$
\lim_{n\rightarrow \omega} |\tau(v_nxv_n^*y)|  \leq \delta < \delta_0, 
$$
for all $x, x_i, y \in F$, $x_0\in F'$. Thus, if we let $v_0=v_n$  for some large enough $n$, then $v_0$ is a unitary element in $P_0'\cap P$ that satisfies 
$\|E_Q(x_0v_0xv_0^*x_0^*)\|_2^2 \leq \delta_0$, for all $x\in F, x_0\in F'$, $|\tau(v_0x_1v_0^*x_2v_0x_3v_0^*x_4)| \leq \delta_0$, for all 
$x_i\in F$, and  $|\tau(v_0xv_0^*y)| \leq \delta_0$, for all $x, y \in F$.

\hfill 
$\square$

\heading 4.  Coarse embeddings of $R$ \endheading

 In this section we prove the main results of the paper, stated in the introduction: Theorem 1.1 (see Theorem 4.2 bellow), Corollaries 
 1.2 and 1.3 (see Corollary 4.3).  

\proclaim{4.1. Lemma} Let $M$ be a $\text{\rm II}_1$ factor, 
$Q\subset M$ a von Neumann subalgebra,  $P \subset M$ an irreducible subfactor such that $P\not\prec_M Q$ and 
$P_0\subset P$ a finite dimensional subfactor. Given  any finite sets $F_0\subset M\ominus P_0$, $F_0'\subset M$, any $m_1\geq 1$ and $\alpha>0$, 
there exists a subfactor $P_1\simeq \Bbb M_{m_1}(\Bbb C)$ in $P_0'\cap P$ such that for any von Neumann subalgebra $P_0^0\subset P_0$ 
we have 
$$
\text{\bf c}_q(y(P^0_0 \vee P_1)y^*, Q) \leq \text{\bf c}_q(yP^0_0y^*, Q) + \alpha,  \forall 1 \leq q \leq 2, y\in F'_0;  \tag a 
$$
$$
|\tau(b_1x_2b_3x_4)|\leq \alpha, \forall b_1, b_3 \in (P_1\ominus \Bbb C1)_1, \forall x_2, x_4 \in F_0; \tag b
$$ 
$$
|\tau(x_0bx_1)| \leq \alpha, \forall b \in (P_1\ominus \Bbb C1)_1, \forall x_0, x_1  \in F'_0. \tag c
$$
\endproclaim
\noindent
{\it Proof}.  Note first that, since $\alpha$ can be taken arbitrarily small independently of the $\| \ \|$-size of the elements in $F_0, F_0'$, 
it is sufficient to prove the statement in the case $F_0$, $F_0'$ are contained in the until ball of $M$, an assumption that we will thus make 
for the rest of the proof. Also, by replacing $F_0$ by $F_0\cup \{x_1x_0-E_{P_0}(x_1x_0) \mid x_0, x_1 \in F_0'\}$, it follows that it is sufficient to 
prove that instead of $(c)$ above, we are required to have 
$$
|\tau(bx)| \leq \alpha, \forall b \in (P_1\ominus \Bbb C1)_1, \forall x \in F_0. \tag c'
$$

We will assume the subalgebra $P_0^0\subset P_0$ is equal to $P_0$, but the proof in this case will in fact show 
that the estimates in $(a)$ hold for any subalgebra $P_0^0\subset P_0$.

Let $1 \in \Cal U_0\subset \Cal U(P_0)$ be an orthonormal basis of $L^2P_0$ made up of unitary elements. Thus, $|\Cal U_0|=m_0^2$, where 
$P_0 \simeq \Bbb M_{m_0}(\Bbb C)$. 

Let $\delta = \alpha^2/36m_0^3m_1$. Let  $P_1^0\subset P_0'\cap P$ be a type I$_{m_1}$ subfactor and $F_1\subset (L^2(P_1^0) \ominus \Bbb C)_1$ be a finite subset  
that's $\delta$ $\| \ \|_2$-dense in $(L^2(P_1^0) \ominus \Bbb C)_1$, i.e., any $x\in (L^2(P^0_1)\ominus \Bbb C)_1$ is $\leq \delta$-close to $F_1$ in the norm $\| \ \|_2$. 

Note that $F := \Cal U_0 (F_0 \cup F_1) \Cal U_0 $ is orthogonal to $P_0$. By applying Lemma 3.1 to this finite set $F$, to $F'=F_0'$ and $\delta_0=\delta > 0$, 
we get a unitary element $v\in P_0'\cap P$ such that 
$$
\|E_Q(yv x v^*y^*)\|_2 \leq \sqrt{\delta}, \forall x\in F, y\in F'\tag 1a 
$$
$$
|\tau(vx_1v^*x_2vx_3v^*x_4)| \leq \delta, \forall x_i\in F \tag 1b
$$
$$
|\tau(vxv^*y)| \leq \delta, \forall x, y \in F. \tag 1c
$$
By taking into account that $F_1$ is $\delta$ $\| \ \|_2$-dense in $(L^2(P_1^0)\ominus \Bbb C)_1$ 
and applying the triangle inequality in $(1a)$ for $x\in F_1\Cal U_0$, it follows that 
$$
\|E_Q(yv x v^*y^*)\|_2 \leq 3\sqrt{\delta}, \forall x\in (L^2(P^0_1) \ominus \Bbb C)_1 \Cal U_0, y\in F'. \tag 2a
$$

Denote $P_1=vP_1^0v^*$ and note that $v((P^0_1 \ominus \Bbb C) \Cal U_0)v^* =(P_1\ominus \Bbb C)\Cal U_0$. 
Thus, by applying (2a) to $x=b u$ where $b\in    (L^2(P_1) \ominus \Bbb C)_1$ and $u\in \Cal U_0$, it follows that  
$$
\|E_Q(ybuy^*)\|_2 \leq 3\sqrt{\delta}, \forall y\in (L^2(P_1) \ominus \Bbb C)_1, u\in  \Cal U_0, y\in F'.
$$
or equivalently 

$$
\sup \{|\tau(ybuy^* z)| \mid z\in (L^2Q)_1\} \leq 3\sqrt{\delta}, \forall b\in  (L^2(P_1) \ominus \Bbb C)_1, u\in  \Cal U_0, y\in F'.  \tag 3a 
$$

Since $\Cal U_0$ is an orthonormal basis for $L^2(P_0)$ and $1\in \Cal U_0$, any element $x\in P_0 \vee P_1$ 
can be uniquely written as $\sum_{u\in \Cal U_0} ub_u$, for some $b_u\in P_1$, with $\|x\|_2^2 = \sum_u \|b_u\|_2^2$. 
Moreover, one has  $E_{P_0}(x)=\sum_u \tau(b_u) u$ with $\tau(x)=\tau(b_1)$. 

Thus, if $x'=\sum_u ub_u$ lies in the unit ball of $L^2(P_0\vee P_1)$ and has $0$ expectation onto $P_0$, then $\tau(y_u)=0$, $\forall u$. 
Moreover, if we denote as usual $p=\frac{q}{q-1}\leq 2$ and take into account that the unit ball of $L^pQ$ is contained in the unit ball of $L^2Q$, then   
from (3a) and the remarks in 2.4 we get the estimates 
$$
\|E_Q(yx'y^*)\|_q = \sup \{|\tau( yx'y^* z)| \mid z\in (L^pQ)_1\}  \tag  4a
$$
$$
\leq \sup \{|\tau(y x' y^*z)| \mid z\in (L^2Q)_1\}
$$
$$
\leq \sum_{u\in \Cal U_0} \sup \{|\tau(yub_uy^* z)| \mid z\in (L^2Q)_1\} 
$$
$$
\leq (\sum_{u\in \Cal U_0} \|b_u\|_2) 3\sqrt{\delta} \leq (\sum_{u\in \Cal U_0} \|b_u\|^2_2)^{1/2} |\Cal U_0|^{1/2} 3\sqrt{\delta}   
$$
$$
= \|x'\|_2 m_0 3\sqrt{\delta}   \leq (m_0m_1)^{1/2}\|x'\|_q m_0 3\sqrt{\delta}, 
$$
where for the last inequalities we have used  the Cauchy-Schwarz inequality 
and the fact that, due to Lemma 2.5.2$^\circ$, we have: 

$$
\sum_u \|b_u\|_2^2=\|x'\|_2^2 \leq \|x'\|_q (\text{\rm dim}(P_0\vee P_1))^{1/4}=(m_0m_1)^{1/2}\|x'\|_q.
$$

By writing any $x\in L^2(P_0\vee P_1)\ominus \Bbb C$ as a sum between its projection on $P_0$ and respectively 
on $(P_1\vee P_0)\ominus P_0$, i.e., $x=(x-E_{P_0}(x)) + E_{P_0}(x)$, and taking into account that $\|E_{P_0}(x)\|_q\leq \|x\|_q$ 
and $\|x-E_{P_0}(x)\|_q \leq 2\|x\|_q$, by applying (4a) to $x'=x-E_{P_0}(x)$ it follows that if such $x$ satisfies $\|x\|_q\leq 1$ then 
$$
\|E_Q(yxy^*)\|_q \leq  \|E_Q(yE_{P_0}(x)y^*)\|_q +  \|E_Q(yx'y^*)\|_q\leq  \text{\bf c}_q(yP_0y^*, Q) + \alpha. \tag 5a
$$

In turn, if we apply $(1b)$ (resp. $(1c)$) to $x_1, x_3 \in F_1$ (resp. $x\in F_1$) and $x_2, x_4 \in F_0 \subset (M\ominus P_0)_1$ (resp. 
$y\in F_0$) and taking into account that $vF_1v^*$ is $\delta$ 
$\| \ \|_2$-dense in $(L^2(P_1) \ominus \Bbb C)_1$, by the triangle inequality (+ Cauchy-Schwarz) we  obtain 
$$
|\tau(b_1x_2b_3x_4)| \leq 3\delta \leq \alpha, \forall b_1, b_3 \in (L^2(P_1) \ominus \Bbb C)_1, x_2, x_4\in F_0. \tag 2b
$$
$$
|\tau(bx)| \leq \alpha, \forall b \in (L^2(P_1) \ominus \Bbb C)_1, x \in F_0. \tag 2c
$$

But $(5a), (2b), (2c)$ are just  conditions $(a), (b), (c')$ required in the statement. 

\hfill 
$\square$

\proclaim{4.2. Theorem} Let $M$ be an arbitrary separable  $\text{\rm II}_1$ factor. Let  $P\subset M$ 
be an irreducible subfactor and $Q\subset M$ a von Neumann subalgebra such that $P\not\prec_M Q$. Given any $\varepsilon >0$, 
$P$ contains a copy of the hyperfinite factor $R\subset P$ 
that's coarse in $M$ $($and thus strongly malnormal and mixing in $M)$, coarse 
with respect to $Q$, and satisfies  
$\text{\bf c}_q(R, Q) \leq \varepsilon$,  $\text{\bf c}_p(Q, R) \leq \varepsilon$, $\forall 1\leq q \leq 2\leq p\leq \infty.$
Moreover, there exists a decreasing sequence of factors $R_n\subset R$ such that $R_n$ is coarse in $M$, $\forall n$, and $\cap_n R_n=\Bbb C1$. 

\endproclaim
\noindent
{\it Proof}. Let  $\{y_k\}_{k\geq 0}, y_0=1,$ be a sequence of elements that's $\| \ \|_2$-dense in $(M)_1$. 

We will construct $R$ as the weak closure 
of the algebra $B$ obtained as union of an increasing sequence of matrix algebras 
$B_m =(\Bbb M_{2}(\Bbb C))^{\otimes m}\simeq  \Bbb M_{2^m}(\Bbb C)$. 
The algebras $B_m$ will be constructed recursively, so that at each step $m$ certain inequalities are satisfied. 

To do this, we need some notations. Thus, we view $\Bbb M_{2}(\Bbb C)$ as  spanned by an orthonormal system of 
selfadjoint unitary elements $\Cal U_0=\{u_0, u_1, u_2, u_3\}$ satisfying $u_0=1$, $u_1u_2=-u_2u_1$ and $u_3=u_1u_2$. 
We denote $J_m=\{0, 1, 2, 3\}^m$ the set of $m$-tuples $j=(j_i)^m_{i= 1}$, with entries  $j_i\in \{0, 1, 2, 3\}$,  
which we view as a subset of the set $J=\{0, 1, 2, 3\}^{(\Bbb N)}$ of infinite sequences $j=(j_i)_{i\geq 1}$ with all but finitely many $j_i$ equal to $0$ 
(where $j=(j_i)^m_{i= 1}\in J_m$ is viewed in $J$ by completing its $i>m$ coordinates with $0$). 
Let $\{u_{i, j_i}\}_{0\leq j_i\leq 3}$ be a copy of $\Cal U_0$ inside $B_{i-1}'\cap B_i\simeq \Bbb M_{2}(\Bbb C)$. 
Let $\Cal U_J\subset B$ be the set of unitary elements $\{u_j\}_{j\in J}$ with $u_j=\Pi_{i\geq 1} u_{i,j_i}$.

We construct $B_m\subset P$ recursively so that $B_0=\Bbb C$ and so that if one denotes 
$F_m=\Cal U_{J_{m-1}}\{y_k - E_{B_{m-1}}(y_k) \mid 1\leq k \leq m\} \ \Cal U_{J_{m-1}}$, $m\geq 1$, then 
the following conditions are satisfied: 
$$
\text{\bf c}_q(y_k(B_j'\cap B_m)y_k^*, Q) \leq \varepsilon (2^{-j-1}-2^{-m-1}), \ 1\leq q \leq 2, \ 0\leq k \leq j \leq m; \tag a
$$
$$
|\tau(b_1x_2b_3x_4)|\leq 2^{-6m-1},   \  b_1, b_3 \in (B_{m-1}'\cap B_m)_1, \  x_2, x_4 \in F_m \vee F_m^*. \tag b
$$
$$
|\tau(bx)| \leq  2^{-6m-1},   \ b \in (B_{m-1}'\cap B_m)_1, \  x \in F_m \vee F_m^*. \tag c
$$

Assume we have constructed the algebras $B_m$ up to $m=n$. By applying Lemma 4.1 to $P_0=B_n$, $m_1=2$, 
$F=F_{n+1} \cup F_{n+1}^*$, $F'=\{y_k\mid 0\leq k \leq n\}$ and $\alpha=\varepsilon 2^{-6n-3}$, 
we get a subfactor $P_1\simeq \Bbb M_{2}(\Bbb C)$ inside $B_n'\cap P$ such that for all $1\leq q \leq 2$, 
$0\leq k \leq j \leq n$  we have 
$$
\text{\bf c}_q(y_k(B_j'\cap B_n  \vee P_1)y_k^*, Q) \leq \text{\bf c}_q(y_k(B_j'\cap B_n)y_k^*, Q)+ \varepsilon 2^{-n-2} 
$$
$$
\leq \varepsilon (2^{-j-1}-2^{-n-1}) + \varepsilon 2^{-n-2}  =  \varepsilon (2^{-j-1}-2^{-n-2}), 
$$
while for all  $ b, b_1, b_3 \in (P_1)_1$ and all  $y, y_2, y_4 \in F_{n+1} \cup F_{n+1}^*$, we have 
$$
|\tau(b_1y_2b_3y_4)|\leq 2^{-6n-3}. 
$$ 
$$
|\tau(by)|\leq 2^{-6n-3}. 
$$

Thus, if we let $B_{n+1}=B_n \vee P_1$, then all conditions $(a), (b), (c)$  are satisfied for $m=n+1$. 

Define $R=\overline{\cup_n B_n}^w$. Then $R$ is a copy of the hyperfinite II$_1$ factor inside $P$, which by (a) and Lemma 2.4.2.3$^\circ$   
satisfies $\text{\bf c}_q(R, Q) \leq \varepsilon/2$, $\forall 1\leq q \leq 2$. By Lemma 2.4.2.1$^\circ$, this implies $\text{\bf c}_p(Q, R) \leq \varepsilon$, $\forall 2\leq p \leq \infty$, 
showing that $R$ satisfies the last condition  of the statement. Also, by Lemma 2.7.1, it follows that the von Neumann algebra generated 
by $R$ and $Q^{op}=J_MQJ_M$ on $L^2M$ is $R\overline{\otimes} Q^{op}$, i.e., $R, Q\subset M$ is a coarse pair. 

Let us prove that condition $(c)$ above implies the following: 

\vskip.05in

{\it Fact}. $\|E_R(y_n)-E_{B_k}(y_n)\|_2 \leq  2^{-3k}$, for all $k\geq n$.

\vskip.05in

To see this, note that for each $m$, the set of unitaries $\{u_j \mid j\in J_m\}$ for an orthonormal basis of $B_m$ while 
$\{u_j \mid j\in J\}$ is an orthonormal basis of $R$. Thus 
$$
E_R(y_n)-E_{B_k}(y_n)=\sum_{j\in J\setminus J_k} \tau(y_nu_j^*)u_j = \sum_{m\geq k} \  \sum_{j\in J_{m+1}\setminus J_m} \tau(y_nu_j^*)u_j 
$$
By Pythagoras Theorem, this gives
$$
\|E_R(y_n)-E_{B_k}(y_n)\|_2^2= \sum_{m\geq k} \sum_{j\in J_{m+1}\setminus J_m} |\tau(y_nu_j^*)|^2 
$$
By property $(c)$,  for each $j\in J_{m+1}\setminus J_m$ and $m\geq k \geq n$ we have $|\tau(y_nu_j^*)| \leq 2^{-6m-1}$. Since $|J_{m+1}\setminus J_m|
=3\cdot 4^{m}$, it follows that 
$$
\|E_R(y_n)-E_{B_k}(y_n)\|_2^2\leq 3\sum_{m\geq k} 2^{-6m-2}= 3\sum_{m\geq 0} 64^{-m}/64^{k}\cdot 4=64^{-k}. 
$$

\vskip.05in

We'll now use $(b)$ and 
the {\it Fact}  above, together with Corollary 2.7.2, to prove that $R$ is coarse in $M$.  
To this end, it is sufficient  to show that for any $n$ and any $\varepsilon >0$, there exists $m$ such that for any $X\in ((B_m'\cap B) \otimes (B_m'\cap B)^{op})_1$ 
with $\tilde{\tau}(X)=0,$ we have 
$$
| \langle X(\xi_n), \xi_n\rangle_{L^2M} | \leq \varepsilon, \tag 1
$$
where $\xi_n=y_n-E_R(y_n)$ is viewed here as a vector  in $L^2M\ominus L^2R\subset L^2M$. 

Let $m$ be so that $\|E_R(y_n)-E_{B_m}(y_n)\|_2 \leq \varepsilon/16$ and $2^{-m}\leq \varepsilon/4$. 

Writing $X\in (B \otimes B^{op})_1$ in the form $X=\underset{j,j'\in J}\to{\sum} c_{j,j'} u_j \otimes u_{j'}$, the condition  $\tilde{\tau}(X)=0$ amounts to 
 $c_{0,0}=0$ and the condition $X\in (B_m'\cap B)\otimes (B_m'\cap B)^{op}$ amounts to $X$ being 
supported by the set $J^m$ of indices $j, j'\in J$ having the first $m$ coordinates equal to $0$. So in order for $(1)$ to be satisfied 
we have to show that  
$$
|\sum_{j,j' \in J^m} c_{j,j'} \tau(u_j \xi_n u_{j'}\xi_n^*)| \leq \varepsilon. \tag 2
$$ 
By the Cauchy-Schwartz inequality, the left hand term is majorized by 
$$
(\sum_{j,j' \in J^m} |c_{j,j'}|^2)^{1/2} (\sum_{j,j' \in J^m} |\tau(u_j \xi_n u_{j'}\xi_n^*)|^2)^{1/2}. 
$$
Since $(\sum_{j,j'} |c_{j,j'}|^2)^{1/2} =\|X\|_{2,\tilde{\tau}}$ in $B\otimes B^{op} \subset R\overline{\otimes} R^{op}$ and since $\|X\|\leq 1$, 
we have $(\sum_{j,j'} |c_{j,j'}|^2)^{1/2} \leq 1$. 
Thus, in order for $(2)$ to be satisfied, it is sufficient to show that: 
$$
\sum_{j,j' \in J^m} |\tau(u_j \xi_n u_{j'}\xi_n^*)|^2 \leq \varepsilon^2 \tag 3
$$

For each $k>m$ we let $J_k^m$ be the set of indices $j\in J$ with $j_i=0$ for all $i>k$ and all $i\leq m$, while for $i=k$ 
one has $j_i\neq 0$. We also let $\overline{J}^m_k=(J^m\cap J_k) \setminus J^m_k$. 
Note that  $J^m$ is the disjoint union of the subsets $J^m_k$ as $k$ runs from $m+1$ to $\infty$, 
with $|J_m^k|=3\cdot 4^{k-m-1}$. We write the sum on the left hand side in $(3)$ as  $\sum_{k\geq m+1} \Sigma(k)$ where 
$$
\Sigma(k) = \sum_{j,j'\in J^m_k}|\tau(u_j\xi_nu_{j'}\xi_n^*)|^2 \tag 4
$$
$$
+\sum_{j\in J^m_k, j'\in \overline{J}^m_k}|\tau(u_j\xi_n u_{j'}\xi_n^*)|^2 +\sum_{j\in \overline{J}^m_k, j'\in J^m_k}|\tau(u_j\xi_n u_{j'}\xi_n^*)|^2. 
$$
Denote by $\Sigma_1(k), \Sigma_2(k), \Sigma_3(k)$ the three sums  on the right hand side above  
and let $y_{n,k}=y_n-E_{B_k}(y_n)$. Since $\xi_n=(y_n-E_{B_k}(y_n))-(E_R(y_n)-E_{B_k}(y_n))=y_{n,k}-(E_R(y_n)-E_{B_k}(y_n)),$ 
for the terms   in $\Sigma_1(k)$ we have the estimates
$$
|\tau(u_j\xi_nu_{j'}\xi_n^*)|\leq |\tau(u_jy_{n,k}u_{j'}y_{n,k}^*)|
$$
$$
+|\tau(u_j(E_R(y_n)-E_{B_k}(y_n))u_{j'}\xi_n^*)|+ |\tau(u_jy_{n,k}u_{j'}(E_R(y_n)-E_{B_k}(y_n))^*)|
$$
$$
\leq  |\tau(u_jy_{n,k}u_{j'}y_{n,k}^*)| + 2\|y_n\|_2\|E_R(y_n)-E_{B_k}(y_n)\|_2, 
$$
where for the last estimate we have used the Cauchy-Schwartz inequality. By condition $(b)$ we have $ |\tau(u_jy_{n,k}u_{j'}y_{n,k}^*)| \leq 2^{-3k}$ 
while by the {\it Fact} above we have $2\|y_n\|_2\|E_R(y_n)-E_{B_k}(y_n)\|_2 \leq 2^{-3k+1}$. 
Since $j,j'\in J^m_k$ and $|J_m^k|=3\cdot 4^{k-m-1}$, it follows that 
$$
\Sigma_1(k) \leq (3\cdot 4^{k-m-1})^2(2^{-3k}+2^{-3k+1})^2 \leq 3^4 2^{-2k-4m-4}\leq 2^{-2k-4m+3} \tag 5
$$

Similarly, by using again that $\xi_n=y_{n,k}-(E_R(y_n)-E_{B_k}(y_n))$, for the terms in $\Sigma_2(k)$ we get the estimate 
$$
|\tau(u_j\xi_n u_{j'}\xi_n^*)| \leq  |\tau(u_jy_{n,k}u_{j'}y_{n,k}^*)| + 2\|y_n\|_2\|E_R(y_n)-E_{B_k}(y_n)\|_2  
$$
$$
\leq   |\tau(u_jy_{n,k}u_{j'}y_{n,k}^*)|+ 2^{-3k+1} \leq 2^{-3k} + 2^{-3k+1}=3\cdot 2^{-3k},  
$$
where we have used that for $j\in J^m_k$ and $j'\in \overline{J}^m_k$ one has $|\tau(u_jy_{n,k}u_{j'}y_{n,k}^*)|\leq 2^{-3k}$ (by condition $(b)$).  
Summing up over $j\in J^m_k$, $j'\in \overline{J}^m_k$ and using the fact that $J^m_k|=3 \cdot 4^{k-m-1}$, $|\overline{J}^m_k|=4^{k-m-1}$, 
it follows that 
$$
\Sigma_2(k)\leq 3^2\cdot 2^{-6k} 3 \cdot 4^{k-m-1} 4^{k-m-1}\leq 3^3 \cdot 2^{-2k-4m-4}\leq 2^{-2k-4m+1}. \tag 6
$$
In exactly the same way, we also get
$$
\Sigma_3(k) \leq  2^{-2k-4m+1}. \tag 7
$$
By adding up $(5), (6), (7)$, we obtain in $(4)$: 

$$
\Sigma(k)\leq 2^{-2k-4m+3}+ 2\cdot 2^{-2k-4m+1}\leq 2^{-2k-4m+4}. \tag 8
$$ 
So after summing up over $k\geq m+1$ 
we get the following estimate for the left hand side of $(3)$: 
$$
(\sum_{j,j' \in J^m} |\tau(u_j \xi_n u_{j'}\xi_n^*)|^2)^{1/2} \leq (\sum_{k\geq m+1} 2^{-2k-4m+4})^{1/2} 
$$
$$
= (\sum_{k\geq 1} 2^{-2k} 2^{-6m+4})^{1/2} \leq 2^{-3m+2} \leq \varepsilon.
$$

To prove the last part of the statement, let us first note that given 
any finite subset of $R$ (which is a coarse hyperfinite subfactor of $M$ that will be fixed from now on, as given by the first part of the proof) and 
any $\varepsilon >0$, there exists a coarse subfactor $R_0\subset R$ such that $F \perp_\delta R_0$. Indeed, since $R$ is AFD, there exists 
a finite dimensional subfactor $N_0\subset R$ such that $F\subset_{\delta/2} N_0$. Take any irreducible subfactor with infinite index $N_1\subset N_0'\cap R$ 
and denote $Q_0=N_0 \vee N_1$, which is thus an irreducible subfactor with infinite index in $R$ so that $F\subset_{\delta/2} Q_0$. By the first 
part of the theorem, there exists a coarse $R_0\subset R$ such that $R_0\perp_{\delta/2} Q_0$. By the triangle inequality, $R_0$ will also 
satisfy $R_0\perp_{\delta} F$. 

We then construct a decreasing sequence of coarse subfactors $R_n$ recursively, as follows. 
Let $\{x_n\}_n\subset (R\ominus \Bbb C1)_1$ be a sequence that's $\| \ \|_2$-dense in the unit ball of $R\ominus1$ and denote $X_n=\{x_1, ..., x_n\}$. 
Assume $R_n\subset R$ is coarse 
and satisfies $X_n \perp_{2^{-n}} R_n$. By the argument above, there exists  a coarse subfactor $R_{n+1} \subset R_n$ 
such that if we let $F=E_{R_n}(X_{n+1})$ then $R_{n+1} \perp_{2^{-n-1}} F$.

Thus, if we denote $D=\cap_n R_n$, then the above conditions imply $\{x_n\}_n \perp D$. Thus $R\ominus \Bbb C1 \perp D$, so that $D=\Bbb C1$. 
\hfill 
$\square$

\proclaim{4.3. Corollary} Let $M$ be a separable $\text{\rm II}_1$ factor,   $P\subset M$ an irreducible subfactor, 
$Q\subset M$ a von Neumann subalgebra satisfying $P\not\prec_M Q$, and $\varepsilon>0$.

$1^\circ$ $P$ contains a MASA $A$ of $M$ that's coarse, mixing and strongly malnormal in $M$, 
has infinite multiplicity, is coarse with respect to $Q$ and satisfies $A\perp_\varepsilon Q$. 

$2^\circ$  There exists a semiregular MASA $D$ of $M$  that's contained in $P$, whose normalizer $\Cal N_M(A)$ lies in $P$ and generates a hyperfinite factor $R$ 
that's coarse with respect to $Q$ and satisfyies $R \perp_\varepsilon Q$.  

$3^\circ$ If $G$ is a countable amenable group, then there exists a copy $\{u_g\}_{g\in G}\subset P$ of the left regular representation of $G$ such that 
$\|E_Q(u_g)\|_2 \leq \varepsilon$, $\forall g\in G\setminus \{e\}$.
\endproclaim 
\noindent
{\it Proof.} $1^\circ$ This is immediate by Theorem 4.2 and Proposition 2.6.3(b) (combined with Example 2.6.4.3$^\circ$). Let us also give a 
self contained  argument. Thus, by Theorem 4.2, $P$ contains a hyperfinite II$_1$ factor $R\subset P$ that's coarse in $M$ and satisfies $R\perp_\varepsilon Q$. 
If $A\subset R$ is any MASA, then $A$ is also a MASA in $M$. Moreover, since $A$ is a subalgebra of $R$, $A \vee A^{op}\subset \Cal B(L^2M\ominus L^2R)$ 
gives a representation of $A\overline{\otimes} A^{op}$ in $\Cal B(L^2(M\ominus R))$, whose commutant contains the commutant of $R\vee R^{op}$ 
which is of type II. Hence,  $(A\overline{\otimes} A^{op})' \cap \Cal B(L^2(M\ominus R))$ is of type I$_\infty$ and therefore, if $A$ is taken to be coarse, 
strongly malnormal, mixing, and with Pukanszky invariant equal to $\infty$ in $R$, like 
in Example 2.6.4.3$^\circ$, then the inclusion $A\subset M$ is coarse, strongly malnormal, mixing, and with Pukanszky invariant equal to $\infty$ as well.

$2^\circ$ Take again $R\subset P$ to be coarse in $M$ and to satisfy $R\perp_\varepsilon Q$. If $D\subset R$ is a Cartan subalgebra, 
then its normalizer in $M$ is contained in $R$, and thus $\Cal N_M(A)''=R\subset P$.

$3^\circ$ By [C76], the left regular representation $\{u_g\}_g$ of any countable amenable group $G$ lies in $R$. So if one embeds $R$ in $P$ 
such that $R\perp_\varepsilon Q$, as in 
$1^\circ$ or $2^\circ$ above, then $\{u_g\}_g$ will satisfy the condition. 
\hfill 
$\square$ 

\proclaim{4.4. Corollary} Let $M$ be a separable $\text{\rm II}_1$ factor,   $P\subset M$ an irreducible subfactor. 

$1^\circ$ $P$ contains a hyperfinite subfactor $R\subset P$ that's coarse with infinite multiplicity  in $M$, i.e., 
it satisfies $_RL^2(M\ominus R)_R\simeq (L^2R \overline{\otimes}L^2R^{op})^{\oplus \infty}$, or equivalently  $R\vee R^{op} \subset \Cal B(L^2(M\ominus R))$ 
is of type $\text{\rm II}_1$ with $\text{\rm II}_\infty$ commutant. 

$2^\circ$ $P$ contains a pair of hyperfinite subfactors $R_0, R_1\subset P$ that's coarse with infinite multiplicity in $M$, i.e. $_{R_0}L^2M_{R_1}\simeq 
(L^2R_0\overline{\otimes} L^2R_1^{op})^{\oplus \infty}$, in other words $R_0\vee R_1^{op} \subset \Cal B(L^2M)$ 
is of type $\text{\rm II}_1$ with $\text{\rm II}_\infty$ commutant. 

\endproclaim
\noindent
{\it Proof.} $1^\circ$ By Theorem 4.2, $P$ contains a copy $S$ of the hyperfinite II$_1$ factor that's coarse in $M$. By Example 2.5.3.3$^\circ$ 
$S$ contains a coarse copy of $R$ with infinite multiplicity, so by Proposition 2.6.3(b), $R$ follows coarse in $M$, and since $(R\vee R^{op})'\cap \Cal B(L^2(M\ominus R))$ 
is a type II factor and has a properly infinite part, it follows that it is II$_\infty$. 

$2^\circ$ Like in the proof of $1^\circ$ above, one first takes a hyperfinite subfactor $S\subset P$ that's coarse in $M$, then let $R_0, R_1\subset S$ be a coarse 
pair in $S$ (cf. Example 2.5.3.4$^\circ$). It follows that $R_0, R_1$ is a coarse pair in $M$. Replacing if necessary $R_0$ by an infinite index subfactor, 
it follows that the type II factor $(R_0\vee R_1^{op})'\cap \Cal B(L^2M)$ is necessarily infinite. 
\hfill 
$\square$ 

\vskip.05in 
\noindent
{\bf 4.5. Remarks.} $1^\circ$ It  remains as an open problem whether for any irreducible inclusion of separable II$_1$ factors  
with infinite index $Q\subset M$ and any $\varepsilon >0$, one can find a hyperfinite II$_1$ factor $R\subset M$ 
such that $\text{\bf c}_p(R, Q)\leq \varepsilon$ for all $1\leq p \leq \infty$ (uniformly in $p$). In particular, whether there exists an irreducible hyperfinite subfactor 
$R\subset M$ such that $\|E_Q(b)\|\leq \varepsilon \|b\|$, for all $b\in R\ominus \Bbb C1$. Note that if true, this would show that for any countable amenable group $G$ 
there exists   a copy $\{u_g\}_{g\in G}$ of the left regular representation of $G$ so that $\|E_Q(u_g)\|\leq \varepsilon$, $\forall  g\in G\setminus \{e\}$.  

$2^\circ$ One can prove various relative versions of 
Theorem 4.2 and Corollary 4.3. For instance, one can show that if $M$ is a II$_1$ factor with a Cartan subalgebra  $A\subset M$, 
then there exists an intermediate hyperfinite subfactor  $A\subset R \subset M$ for which there exist 
unitaries in the normalizer of $A$ in $M$, $\{u_n\}_n \subset \Cal N_M(A)$, such that $L^2M\ominus L^2R=\oplus_n L^2(Ru_nR)$, 
with $L^2(Ru_nR)\simeq L^2R\overline{\otimes} L^2R^{op}$, $\forall n$.  
Along these lines, it would be interesting to  see whether if $A\subset Q\subset M$ is an intermediate subfactor 
with infinite index in $M$, contains a Cartan subalgebra $A$ of $M$, then given any amenable group $\Gamma$,  there exists a copy 
of the left regular representation of $\Gamma$,  
$\{u_g\}_{g\in \Gamma}$, in the normalizer of $A$ in $M$, so that $E_Q(u_g)=0$, $\forall g\in G\setminus \{e\}$. 

$3^\circ$  Note that, by using  rather minimal additional effort, one can prove the following stronger form of Theorem 4.3. 
Let $M$ be a separable $\text{\rm II}_1$ factor and $P\subset M$ an irreducible subfactor. Let $Q_n \subset M$ be a sequence of 
von Neumann subalgebras such that $P\not\prec_M Q_n$, $\forall n$. There exists  
a hyperfinite subfactor $R \subset P$ that's coarse in $M$ and  is also coarse with respect to $Q_n$, for all $n$.  Moreover, if $Q_n$ is a finite family, 
then given any $\varepsilon >0$, one can construct $R\subset P$ so that in addition to the above properties one has 
$R\perp_\varepsilon Q_n$,  $\forall n$. However, if the family $Q_n$ is infinite, then one cannot expect to have this latter 
property uniformly in $\varepsilon$, simultaneously for all $Q_n$. For instance, one can take $P=M$ and $Q_n$ an increasing family 
of irreducible subfactors with infinite index such that $Q_n \nearrow M$. In such a case, the condition $R\perp_{\varepsilon} Q_n$, $\forall n$, 
would imply $R\perp_\varepsilon M$, a contradiction. 

$4^\circ$ Note that the control over the asymptotic 2-independence in the iterative construction of $R$ in the proof of Theorem 4.2   
is first made with respect to the Hilbert norm $\| \ \|_2$, at each ``local level''.  
But in order to get coarseness at the end (which requires representing $R\overline{\otimes} R^{op}$), one has to covert it to $\| \ \|_1-\| \ \|_\infty$ estimates 
(using  Lemma 2.4.2), which is why we need the asymptotic 2-independence to converge ``very  fast''. If we only wanted $R\subset M$ to be  mixing at the end 
(which by Section 2.6 is a slightly weaker property than coarseness) then the $\| \ \|_2$-estimates at local levels are easier. 
One interesting aspect about mixing is that  it has a  natural ``multiple'' version, as defined below. 

\vskip.05in
\noindent
{\bf 4.6. Definition}. Let $B\subset M$ be a diffuse von Neumann subalgebra. 
We say that $\Cal U(B) \curvearrowright^{\text{\rm Ad}} M$ is $n$-{\it mixing relative 
to} $B$ if 
$$
\underset{u_1, ...,u_n \rightarrow 0}\to{\lim}  \|E_B(x_0u_1x_1... u_nx_n)\|_2 = 0,  
\forall x_0, x_1, ..., x_n\in M\ominus B,
$$ 
where the limit is over $u_i\in \Cal U(B)$ tending weakly to $0$. We then also say that 
$B\subset M$ is $n$-{\it mixing}. If this property holds true for all $n\geq 1$, then we say that $B\subset M$ 
has the {\it multiple mixing} property. 

\vskip.05in 

An obvious example of  multiple mixing inclusions $B\subset M$ is when $B$ is freely complemented in $M$, i.e., $M=B*N_0$, with $B$ diffuse. 
In fact, by its very formulation, the multiple mixing property for an inclusion $B\subset M$ can be viewed as $B$ being {\it almost freely complemented} in $M$.

If $\Gamma \curvearrowright (N,\tau)$ is a free 
action, then one can easily show that this action is $n$-mixing iff $L(\Gamma) \subset M=N\rtimes \Gamma$ 
is $n$-mixing (one implication is trivial, by applying condition 4.6 to $u_i=u_{g_i}$ where $g_i\in \Gamma$ are elements of the group 
that go to $\infty$; we leave the other implication as an exercise). Thus, if such an action is $n$-mixing for all $n$, then this inclusion follows 
multiple mixing. Since a Bernoulli $\Gamma$-action with base $(N_0, \tau_0)$ is $n$ mixing, $\forall n$ 
(see e.g. page 472 in [P03]), it follows that the resulting inclusion $L(\Gamma)\subset M=N_0^{\overline{\otimes}\Gamma} \rtimes \Gamma$ is multiple mixing.

The arguments used in the proof of Theorem 4.2 can be easily adapted to prove the following: 

\proclaim{4.7. Theorem} If $M$ is a separable $\text{\rm II}_1$ factor and $P\subset M$ an irreducible subfactor, then $P$ 
contains a copy of the hyperfinite factor $R$  that's multiple mixing in $M$.  

\endproclaim

\head 5. On the coarse nature of free group factors \endhead

The iterative technique for constructing embeddings under constraints of AFD algebras into a given II$_1$ factor $M$ that we have developed 
in this paper suggests that, under appropriate hypothesis on $M$, one should be able to construct pairs of embeddings of the hyperfinite factor, $R_0, R_1 \hookrightarrow M$, 
so that the corresponding bimodule $_{R_0}L^2M_{R_1}$  satisfies certain properties, as a consequence of the  
``local steps'' taken in the recursive procedure of building $R_0, R_1$.

There are two sets of problems that come out quite naturally from this idea, both having to do with the free group factors.

The first problem is about showing that $L(\Bbb F_\infty)$ cannot be generated 
by finitely many elements. Recall from (Corollary 4.7 in [R92]) that this fact would also imply that the free group factors $L(\Bbb F_n)$, $2\leq n \leq \infty$, are all non-isomorphic. 

We will say that a II$_1$ factor $M$ is {\it stably single generated}, or that it has the {\it stable single generation} ({\it SSG}) property, if there exist $t_n \searrow 0$ such that $M^{t_n}$ 
is single generated, $\forall n$ (N.B. this is easily seen to be equivalent to $M^t$ being single generated, $\forall t>0$).

\vskip.05in

\noindent
{\bf 5.1. Conjecture}. $(a)$ {\it If a } II$_1$ {\it factor $M$ has the SSG property then it contains a pair 
of hyperfinite factors $R_0, R_1 \subset M$ such that $R_0\vee R_1^{op}$ is a purely infinite 
von Neumann algebra}.  

$(b)$ {\it If a } II$_1$ {\it factor $M$  has the SSG property, then it contains a pair 
of hyperfinite factors $R_0, R_1\subset M$ such that $R_0 \vee R_1^{op}=\Cal B(L^2M)$}.  

\vskip.05in  

Of course, part $(b)$ above is stronger than $(a)$. While we comment in more details on the motivations behind this problem in 
other papers (see e.g. Section 7 in [P19]), 
let us point out right away that  if $L(\Bbb F_\infty)$ would be finitely generated, then  by using the fact that  it has non-trivial fundamental 
group (by [V89]; cf. also  [R92]), it would be SSG. If $(a)$ above holds true and $R_0, R_1\subset  M=L(\Bbb F_\infty)$ are hyperfinite subfactors so that $R_0 \vee R_1^{op}$ is properly 
infinite in $\Cal B(L^2M)$, then $R_0\vee R_1^{op}$ would have a cyclic vector (see e.g. [D57]), i.e., there would exist $\xi\in L^2M$ 
such that $[R_0 \xi R_1]=L^2M$. But this contradicts (Theorem 4.2 in [GeP98]). Thus, if true, Conjectures 5.1 would imply that $L(\Bbb F_\infty)$   cannot be 
generated by finitely many elements, which in turn would imply by ([R92]) that all free group factors $L(\Bbb F_n), 2\leq n \leq \infty$, are non-isomorphic.

The second set of problems is inspired by the fact that  free group factors $L(\Bbb F_n)$ are particularly prompt to the type of dichotomic decomposition 
into ``coarse plus trivial'' bimodules over amenable subalgebras, 
that we emphasized in this paper. This pattern  has been 
first recognized in ([P81b], [P81d])     
and was much exploited in deformation-rigidity theory (see e.g., [P01], [IPeP05], [OP07], [PV11], [I12], etc). In particular,  
such decompositions were key to establishing amenability properties for normalizers of amenable subalgebras  of $L(\Bbb F_n)$ in ([OP07]). 

Another manifestation of this paradigm is  the striking 
structural result discovered  by Voiculescu, using his free entropy dimension theory: the bimodule decomposition of 
$L(\Bbb F_n)$ over any abelian von Neumann subalgebra $B$ must contain a copy of coarse bimodule $L^2B\overline{\otimes}{L^2B}^{op}$   
(cf. [V96]; see also [GeP98], [Ju07], [Ha15]). All this leads to the following:

\vskip.05in
\noindent
{\bf 5.2. Conjecture}. $(a)$ {\it Any maximal amenable $($equivalently maximal AFD$)$ von Neumann subalgebra $B$ 
of $L(\Bbb F_t)$ is coarse, $\forall 1< t \leq \infty$} $($cf. also 4.3 in $\text{\rm [Ha15]})$. 

$(b)$ {\it If $B, B_0 \subset L(\Bbb F_t)$, $1< t \leq \infty$,  are maximal amenable von Neumann subalgebras, then there exist projections 
$p\in B, p_0\in B_0$ such that $pBp$ and $p_0B_0p_0$ are unitary conjugate in $L(\Bbb F_n)$, while $(1-p)B(1-p), (1-p_0)B_0(1-p_0)$ 
is a coarse pair.}

\vskip.05in

We will refer to 5.2 above as the {\it coarseness conjecture} for free group factors. \footnote{After an early version of this paper has been circulated, we learned from Ben Hayes 
that he conjectured in (4.3 of [Ha15]) that if $B\subset L(\Bbb F_t)$ is maximal amenable, 
then $_BL^2(L(\Bbb F_t))_B \subset L^2(B)^{\oplus \infty}$, which in our terminology amounts to $B$ being coarse in $L(\Bbb F_t)$. 
Thus, part $(a)$ of 5.2 above has already been stated in [Ha15], in an equivalent form.} 
Similarly, the conjecture obtained from 5.2 by replacing 
everywhere ``coarse'' by ``mixing'' (respectively ``strongly malnormal''), will be called the {\it mixingness conjecture} (respectively 
{\it strong malnormality conjecture}). 

A related conjecture formulated in [PeT07] (see very last paragraph in that paper), predicts that diffuse amenable 
von Neumann subalgebras in $L(\Bbb F_t)$ have unique maximal amenable extension. In other words,   
if $B, B_0\subset L(\Bbb F_t)$ are maximal amenable subalgebras, then $B\cap B_0$ diffuse implies $B=B_0$. 

We notice below 
that this is equivalent to a {\it malnormality conjecture}, asserting that any maximal amenable von Neumann subalgebra $B$ in an interpolated free group 
factor $L\Bbb F_t$ is malnormal (i.e., if $u\in \Cal U(L\Bbb F_t)$ satisfies $uBu^*\cap B$ diffuse, then $u\in B$). 
We also show that the strong malnormality, 
mixingness and coarseness give increasingly stronger conjectures.

\proclaim{5.3. Proposition} $1^\circ$ If part $5.2.(a)$ of the coarseness $($resp. mixingness, resp. malnormality$)$ conjecture holds true, then 
part $5.2.(b)$ holds true as well. 

$2^\circ$ Coarseness conjecture $\Rightarrow$ mixingness conjecture $\Rightarrow$ strong malnormality conjecture 
$\Rightarrow$ malnormality conjecture. Moreover, the malnormality conjecture is equivalent to the $\text{\rm [PeT07]}$-conjecture.  
\endproclaim  
\noindent
{\it Proof}. $1^\circ$ Assume any maximal amenable subalgebra of $L(\Bbb F_t)$ is coarse, $\forall 1<t\leq \infty$. Let  
$B, B_0\subset L(\Bbb F_t)$ be maximal amenable and take $\tilde{B}=B\oplus B_0 \subset \Bbb M_{2}(L\Bbb F_t)=L(\Bbb F_{s})$, 
where $s=(t+3)/4$ (cf [R92], [Dy93]). Denote $\{e_{ij}\mid 1\leq i,j \leq 2\}$ the system of matrix units in $\Bbb M_2(\Bbb C)\subset \Bbb M_{2}(L\Bbb F_t)=L(\Bbb F_{s})$. 
Let $\Cal B\subset L(\Bbb F_s)$ be a maximal amenable subalgebra containing $\tilde{B}$. Note that $e_{ii}\in \tilde{B}\subset \Cal B$ and  
$e_{11}\Cal Be_{11}=Be_{11}$, $e_{22}\Cal Be_{22}=B_0e_{22}$. Hence, if $v\in e_{11}\Cal Be_{22}$ 
is a partial isometry with left support $pe_{11} \in Be_{11}$ and right support $p_0e_{22} \in B_0e_{22}$, 
then we necessarily have $v^*(pBpe_11)v=p_0B_0p_0e_{22}$. Moreover, by its form it follows that $v=pue_{12}$ for some unitary $u\in L(\Bbb F_t)$, 
with $u^*pBpu=p_0B_0p_0$. 

Let $(p,p_0)$ be a pair of projections with $p\in B, p_0\in B_0$,  
$p\sim_{\Cal B} p_0$, and such that $(p, p_0)$ is maximal with these properties. From the above it 
follows that $pBp, p_0B_0p_0$ are unitary conjugate in $L(\Bbb F_t)$, 
while by maximality we have $(1-p)e_{11}\Cal Be_{22}(1-p_0)=0$.  But then the coarseness of $\Cal B\subset L(\Bbb F_s)$ 
implies that $L^2((1-p)L^2(\Bbb F_t)(1-p_0))$ is coarse as a $(1-p)B(1-p), (1-p_0)B_0(1-p_0)$  Hilbert bimodule. 

The implication $5.2.(a) \Rightarrow 5.2.(b)$ in the mixing (resp strongly malnormal) case is very similar, using this same $2$-by-$2$ matrix trick. 
We leave details as an exercise. 

\vskip.05in

$2^\circ$  The implications in the first part are immediate by Proposition 2.6.3. To see that  the Peterson-Thom conjecture is equivalent to 
the malnormality conjecture one uses the $2$-by-$2$ matrix trick as above. Indeed, if $B, B_0\subset L(\Bbb F_t)$ are maximal amenable 
and we let $\tilde{B}=B\oplus B_0\subset \Bbb M_{2 \times 2}(L\Bbb F_t)=L(\Bbb F_s)$, 
where $s=(t+3)/4$, then take $\Cal B$ to be a maximal amenable von Neumann subalgebra in $L(\Bbb F_s)$ that contains $\tilde{B}$ and denote 
$u=e_{12}+e_{21}$. Since the condition $B\cap B_0$ diffuse is equivalent to $u\tilde{B}u^* \cap \tilde{B}$ diffuse and $B=B_0$ is equivalent to $e_{12}\in \Cal B$, 
it follows that malnormality of maximal amenable subalgebras is equivalent to the unique maximal amenable extension property in [PeT07].    
\hfill $\square$ 

\vskip.1in 

The techniques in this paper suggest the following approach to the coarseness conjecture. Let  
$B \subset M=L(\Bbb F_n)$ be a maximal amenable von Neumann subalgebra and assume for simplicity $B\simeq R$.  If we assume 
$_B L^2M_B$ has a piece that's not coarse, then that part is not finite. Due to maximal amenability, this should allow to construct 
recursively another copy of $R$ inside $M$, with $_R L^2M_R$ having no coarse component at all (i.e., $R \vee R^{op}$ is purely infinite), 
contradicting ([GeP98]). 

Let us finally mention that in view of Definition 4.5 and Theorem 4.6, it is tempting to 
predict that any maximal amenable  subalgebra $B$ of a free group factor $M=L\Bbb F_n$ is 
multiple mixing (almost freely complemented) in $L\Bbb F_n$. 
More support towards this 
comes from results in [St04].  

\vskip.05in
\noindent
{\bf 5.4. Conjecture}. {\it Any maximal amenable von Neumann subalgebra $B$ of $L(\Bbb F_t)$ is multiple mixing,  $\forall 1< t \leq \infty$}.

\vskip.05in
Conjectures 5.2 and 5.4 are obviously true when the subalgebra $B$ is freely complemented in $M=L\Bbb F_t$ (N.B.: by [P81d], 
if $B$ is amenable diffuse, then $B$ is maximal amenable in $M=B* N_0$, for any tracial von Neumann algebra $N_0 \neq \Bbb C1$). 

Another example of maximal amenable subalgebra in $L\Bbb F_n$, $2\leq n <\infty$, is the so-called
{\it radial} (or Laplacian) MASA, $L_n\subset L\Bbb F_n$, generated by the selfadjoint operator obtained by summing up 
the $n$ free generators  and their inverses (cf. [CFRW09]; see also [R89],  [W15]). 
Calculations in [R89] show that this maximal amenable subalgebra $L_n\subset L\Bbb F_n$ is indeed coarse and multiple mixing. 
The following questions are perhaps ``too bold'', but surprisingly enough, they  are open: \footnote{I am grateful to Dima Shlyakhtenko for pointing out 
to me these questions.}

\vskip.05in
\noindent
{\bf 5.5. Questions}.  {\it Is $L_n$ freely complemented in $L(\Bbb F_n)$, $2\leq n < \infty$} ? {\it Do there exist maximal amenable von Neumann 
subalgebras of $L(\Bbb F_n)$ that are not freely complemented} ?   

\vskip.05in 

If indeed any maximal amenable subalgebra is freely complemented, then this would of course imply both 5.2 and 5.4. To start with, it should be true that 
if $B_0, B_1\subset L(\Bbb F_n)$ are maximal amenable that have no intertwining parts (so they should be mutually coarse, by $5.2(b)$), 
then any direct sum $pB_0p + (1-p)B_1^0(1-p)$ is also freely complemented, where $p\in \Cal P(B_0)$  and  $(1-p)B_1^0(1-p)$ 
is a unitary conjugate of $p_1B_1p_1$, for some $p_1\in \Cal P(B_1)$ of same trace as $1-p$. In other words, direct sum of mutually coarse 
freely complemented diffuse amenable subalgebras in free group factors, should be freely complemented (all this in appropriate 
amplifications). Note in this respect the recent results in [Je19], which show in particular that certain perturbations 
of the free generators in $L(\Bbb F_n)$ are still freely complemented. 

An interesting related question is to classify inclusions of the form $B_0 \subset B_0* N \simeq L(\Bbb F_2)$, where 
$(B_0,\tau_0)$ is a fixed diffuse amenable tracial von Neumann algebra (for instance, $B_0$  taken to be  the diffuse abelian von Neumann algebra $A=L^\infty([0,1])$) 
and $(N, \tau)$ runs over the diffuse amenable tracial von Neumann algebras (note that the isomorphism $B_0* N \simeq L(\Bbb F_2)$ is due to [Dy93]). 
Are these inclusions classified by the isomorphism class of $(N, \tau)$? 
For instance, are the inclusions $A \subset A * A$  and $A \subset A * R$ non-isomorphic?

\head  References \endhead

\item{[AP17]} C. Anantharaman, S. Popa: ``An introduction to II$_1$ factors'', \newline www.math.ucla.edu/$\sim$popa/Books/

\item{[CFRW09]} J. Cameron, J. Fang, M. Ravichandran, S. White:  
{\it The radial masa in a free group factor is maximal injective}, J. London Math. Soc.,  {\bf 72} (2010), 787-809.

\item{[C76]} A. Connes: {\it Classification of injective factors}, Ann. of Math., {\bf 104} (1976), 73-115.

\item{[D54]} J. Dixmier: {\it Sousanneaux abeliens maximaux dans les facteurs de type fini}, Ann. of Math. {\bf 59} (1954), 279-286. 

\item{[D57]} J. Dixmier:  ``Les alg\'ebres d'op\'erateurs dans l'espace Hilbertien,'' Gauthier-Villars, Paris, 1957.

\item{[Dy92]} K. Dykema: {\it Interpolated free group factors}, Pac. J. Math. {\bf 163} (1994), 123-135.

\item{[Dy93]} K. Dykema: {\it Free products of hyperfinite von Neumann algebras and free dimension}, 
Duke Math. J. {\bf 69} (1993), 97-119.

\item{[GP14]} A. Galatan, S. Popa: {\it Smooth bimodules and cohomology of} II$_1$ {\it factors}, Journal of Inst. Math. Jussieu 
(2015), 1-33 (math.OA/1406.6182)

\item{[GeP98]} L. Ge, S. Popa: {\it On some decomposition properties for factors of type} II$_1$,  Duke Math. J., {\bf 94} (1998), 79-101.

\item{[H15]} B. Hayes:  1-{\it bounded entropy and regularity problems in von Neumann algebras},   Int. Math. 
Res. Notices, {\bf 3} (2018), 57-137 (arXiv:1505.06682). 

\item{[I12]} A. Ioana:  {\it Cartan subalgebras of amalgamated free product} II$_1$ {\it factors}. 
(Appendix by A. Ioana and  S. Vaes),  Ann. Sci. Ec. Norm. Super. {\bf 48} (2015),  71-130.

\item{[IPeP05]} A. Ioana, J. Peterson, S. Popa:
{\it Amalgamated Free Products of w-Rigid Factors and Calculation of their
Symmetry Groups}, Acta Math. {\bf 200} (2008), No. 1, 85-153. (math.OA/0505589)

\item{[Je19]} D. Jekel: {\it Conditional Expectation, Entropy, and Transport for Convex Gibbs Laws in Free Probability}, 
 arXiv:1906.10051 

\item{[J81]} V.F.R. Jones: {\it A converse to Ocneanu's theorem}, J. Operator Theory {\bf 10} (1983), 61-63. 

\item{[J83]} V.F.R. Jones: {\it Index for subfactors}, Invent. Math. {\bf 72} (1983), 1-25. 

\item{[Ju07]} K. Jung: {\it Strongly} 1-{\it bounded von Neumann algebras}. Geom. Funct. Anal., {\bf 17} (2007), 1180-1200. 

\item{[MvN43]} F. Murray, J. von Neumann: {\it Rings of operators
IV}, Ann. Math. {\bf 44} (1943), 716-808.

\item{[OP07]} N. Ozawa, S. Popa: {\it On a class of} II$_1$ {\it
factors with at most one Cartan subalgebra}, Annals of Mathematics {\bf 172} (2010),
101-137 (math.OA/0706.3623)

\item{[PeT07]} J. Peterson, A. Thom: {\it Group cocycles and the ring of affiliated operators}, Invent. Math. {\bf 185} (2011), 561-592.

\item{[PiP84]} M. Pimsner, S. Popa: {\it Entropy and index for
subfactors}, Annales Scient. Ecole Norm. Sup., {\bf 19} (1986),
57-106.

\item{[P81a]} S. Popa: {\it On a problem of R.V. Kadison on maximal
abelian *-subalgebras in factors}, Invent. Math., {\bf 65} (1981),
269-281.

\item{[P81b]} S. Popa: {\it Orthogonal pairs of *-subalgebras in
finite von Neumann algebras}, J. Operator Theory, {\bf 9} (1983),
253-268.

\item{[P81c]} S. Popa: {\it Singular maximal abelian *-subalgebras in
continuous von Neumann algebras}, J. Funct. Analysis, {\bf 50}
(1983), 151-166. 

\item{[P81d]} S. Popa: {\it Maximal injective subalgebras in factors
associated with free groups}, Advances in Math., {\bf 50} (1983), 27-48.

\item{[P90]} S. Popa: {\it Markov traces on Universal Jones algebras
and subfactors of finite index}, Invent.  Math., {\bf 111} (1993), 375-405. 

\item{[P92]} S. Popa: {\it Free independent sequences in type} II$_1$ {\it factors
and related problems}, Asterisque, {\bf 232} (1995), 187-202.

\item{[P01]} S. Popa: {\it Some rigidity results for
non-commutative Bernoulli shifts}, J. Fnal. Analysis {\bf 230}
(2006), 273-328 (MSRI preprint No. 2001-005)

\item{[P03]} S. Popa: {\it Strong Rigidity of }  II$_1$ {\it Factors
Arising from Malleable Actions of $w$-Rigid Groups} I, II, Invent. Math.,
{\bf 165} (2006), 369-408 and 409-450.

\item{[P04]} S. Popa: {\it Some computations of $1$-cohomology groups
and construction of non orbit equivalent actions}, Journal of the
Inst. of Math. Jussieu {\bf 5} (2006), 309-332 (math.OA/0407199).

\item{[P05a]} S. Popa: ``Deformation-rigidity theory'', NCGOA mini-course, Vanderbilt University, May 2005.  

\item{[P05b]} S. Popa: {\it Cocycle and orbit equivalence superrigidity for malleable actions of w-rigid groups}, Invent. Math. 
{\bf 170} (2007), 243-295 (math.GR/0512646).

\item{[P06]} S. Popa: {\it Deformation and rigidity for group actions
and von Neumann algebras}, in ``Proceedings of the International
Congress of Mathematicians'' (Madrid 2006), Volume I, EMS Publishing House,
Zurich 2006/2007, pp. 445-479.

\item{[P13a]} S. Popa: {\it A} II$_1$ {\it factor approach to the Kadison-Singer problem}, 
Comm. Math. Physics. {\bf 332} (2014), 379-414 (math.OA/1303.1424).

\item{[P13b]}  S. Popa: {\it Independence properties in subalgebras of ultraproduct} II$_1$ {\it factors}, Journal of Functional Analysis 
{\bf 266} (2014), 5818-5846 (math.OA/1308.3982)

\item{[P16]}  S. Popa:  {\it Constructing MASAs with prescribed properties}, Kyoto J. of Math, {\bf 59} (2019), 367-397 (math.OA/1610.08945)

\item{[P17]} S. Popa:  {\it Asymptotic orthogonalization of subalgebras in type} II$_1$ {\it factors}, Publ. RIMS, Kyoto Univ. {\bf 55} (2019),  795-809 
(math.OA/1707.07317) 

\item{[P18]} S. Popa:  {\it On the vanishing cohomology problem for cocycle actions of groups on} II$_1$ {\it factors},  to appear in Ann. Ec. Norm. Sup., math.OA/1802.09964

\item{[P19]} S. Popa: {\it Tight decomposition of factors and the single generation problem},  math.OA/1910.14653.

\item{[PV11]} S. Popa, S. Vaes: {\it Unique Cartan decomposition for} II$_1$ 
{\it factors arising from arbitrary actions of free groups}, Acta Mathematica, {\bf 194} (2014), 237-284.  

\item{[Po67]} R. Powers: {\it Representations of uniformly hyperfinite algebras and their associated von Neumann rings}, Ann. Math. {\bf 86} (1967), 138-171.

\item{[Pu60]} L. Pukanszky: {\it On maximal abelian subrings in factors of type} II$_1$,  Canad. J. Math. {\bf 12} (1960), 289-296. 

\item{[R89]} F. Radulescu: {\it Singularity of the radial subalgebra of} $L (\Bbb F_N)$ {\it and the Pukanszky invariant},  
Pacific J. Math., {\bf 151} (1991), 297Ð306.

\item{[R92]} F. Radulescu:  {\it Random matrices, amalgamated free products and subfactors of the von Neumann algebra of a free group, of noninteger 
index}, Invent. Math. {\bf 115} (1994), 347-389.

\item{[St04]} M. Stefan: {\it Indecomposability of the free group factors over nonprime subfactors and abelian subalgebras}, 
Pac. J. Math. {\bf 219} (2005), 365-390. 

\item{[V89]} D. Voiculescu: {\it Circular and semicircular systems and free product factors}, 
In ``Operator algebras, unitary representations, enveloping algebras, and invariant theory'' (Paris, 1989), 
Progr. Math. {\bf 92}, Birkh\"{a}user, Boston, 1990, pp. 45-60.

\item{[V96]} D. Voiculescu: {\it The analogues of entropy and of FisherÕs information measure in free probability theory} III: 
{\it The absence of Cartan subalgebras}, Geom. Funct. Anal. {\bf 6} (1996) 172-199.

\item{[W15]} C. Wen:  {\it Maximal amenability and disjointness for the radial masa}, J. Fnal. Analysis {\bf 270} (2016), 787-801. 

\enddocument